\documentclass[12pt]{amsart}
\usepackage{amsmath,amsfonts,amsthm,amscd}
\usepackage{graphicx}
\usepackage{esint}
\usepackage{color}
\def\barint{\kern4pt
\raise3.4pt\hbox{\vrule height.8pt width5pt}%
\kern-9pt % -(4pt + 5pt)
\int}
\newtheorem{theorem}{Theorem}

\newtheorem{proposition}[theorem]{Proposition}

\newtheorem{corollary}[theorem]{Corollary}
\newtheorem{remark}[theorem]{Remark}

\renewcommand{\a }{\alpha }
\renewcommand{\b }{\beta }
\renewcommand{\d}{\delta }
\newcommand{\D }{\Delta }

\newcommand{\e }{\varepsilon }

\newcommand{\G }{\Gamma }

\newcommand{\n }{\nabla }

\newcommand{\s }{\sigma }

\renewcommand{\o }{\omega }

\newcommand{\pa }{\partial}

\newcommand{\ov}{\overline}

\newcommand{\be}{\begin{equation}}
\newcommand{\ee}{\end{equation}}
\newenvironment{pf}{\noindent{\sc Proof}.\enspace}{\rule{2mm}{2mm}\medskip}

\numberwithin{equation}{section} \numberwithin{theorem}{section}
\begin{document}

\newtheorem{lem}{Lemma}[section]
\newtheorem{pro}[lem]{Proposition}
\newtheorem{thm}[lem]{Theorem}
\newtheorem{rem}[lem]{Remark}
\newtheorem{cor}[lem]{Corollary}
\newtheorem{df}[lem]{Definition}

\bibliographystyle{amsalpha}
\title [A strong maximum principle]{A strong maximum principle for the Paneitz operator and a non-local flow for the $Q$-curvature}

\author{Matthew J. Gursky}
\address{Department of Mathematics \\
University of Notre Dame \\
255 Hurley Hall \\
 Notre Dame, IN 46556}
\email{mgursky@nd.edu}
\author{Andrea Malchiodi}
\address{Area of MAthematics \\
         SISSA  \\
         Via Bonomea 265  \\
         34136 Trieste \\
         ITALY \\ and University of Warwick \\ Mathematics Institute - Zeeman Building \\
         Coventry CV4 7AL }
\email{malchiod@sissa.it, A.Malchiodi@warwick.ac.uk}

\begin{abstract}
In this paper we consider Riemannian manifolds $(M^n,g)$ of dimension $n \geq 5$, with semi-positive $Q$-curvature and non-negative scalar curvature.
Under these assumptions we prove $(i)$ the Paneitz operator satisfies a strong maximum principle; $(ii)$ the Paneitz operator is a positive operator;
and $(iii)$ its Green's function is strictly positive.  We then introduce a non-local flow whose stationary points are metrics of constant positive
$Q$-curvature.  Modifying the test function construction of Esposito-Robert, we show that it is possible to choose an initial conformal metric so that
the flow has a sequential limit which is smooth and positive, and defines a conformal metric of constant positive $Q$-curvature.
\end{abstract}

\maketitle
%\date{November, 2013}

%%%%%%%%%%%%%
\section{Introduction}
%%%%%%%%%%%%%%%%

In 1983 S. Paneitz introduced a fourth-order conformally invariant differential operator acting on smooth functions, which is defined on any pseudo-Riemannian manifold \cite{Pan}.  Subsequently, T. Branson \cite{BranSM} recognized that this operator describes the conformal transformation of a curvature quantity which is fourth order in the metric.

To describe the operator and associated curvature quantity, let $A$ denote the Schouten tensor
\begin{align} \label{Adef}
A = \frac{1}{(n-2)}\big( Ric - \frac{1}{2(n-1)} R g \big),
\end{align}
where $Ric$ is the Ricci tensor and $R$ the scalar curvature, and $\sigma_k(A)$ denote the $k^{th}$-symmetric function of the eigenvalues of $A$.  Then the {\em $Q$-curvature} of Branson
is defined by
\begin{align} \label{Qdef}
Q = - \Delta \sigma_1(A) + 4 \sigma_2(A) + \frac{n-4}{2} \sigma_1(A)^2,
\end{align}
and the eponymous operator of Paneitz is
\begin{align} \label{Pdef}
P_g u = \Delta_g^2 u + \mbox{div}_g \big\{ \big( 4 A_g - (n-2) \sigma_1(A_g)g \big)(\nabla u, \cdot) \big\} + \frac{n-4}{2} Q_g u.
\end{align}
The formula connecting $P$ to $Q$ is the following: if the dimension $n \neq 4$, suppose $\hat{g} = u^{\frac{4}{n-4}}g$ is a conformal metric; then the $Q$-curvature of $\hat{g}$ is given by
\begin{align} \label{PQ}
Q_{\hat{g}} = \frac{2}{n-4} u^{-\frac{n+4}{n-4}} P_g u.
\end{align}
When the dimension is four one writes $\hat{g} = e^{2w}g$, and
\begin{align} \label{P4Q}
Q_{\hat{g}} = e^{-4w} \big( - \frac{1}{2} P_g w + Q_g \big).
\end{align}
%Curiously, Branson also pointed out that (\ref{P4Q}) is the limiting expression of (\ref{PQ}), if one formally views the dimension as a real number.

Branson pointed out that the formulas (\ref{PQ}) -- (\ref{P4Q}) naturally suggest a higher order version of the Yamabe problem: given $(M^n,g)$, find a conformal metric of constant $Q$-curvature. In dimensions $n \neq 4$ this is equivalent to finding a \underline{positive} solution of
\begin{align} \label{BPE}
P_g u = \lambda \, u^{\frac{n+4}{n-4}},
\end{align}
where $\lambda$ is a constant. In four dimensions the equation is
\begin{align} \label{PB4}
P_g w + 2 Q_g = \lambda \, e^{4w}.
\end{align}
In both cases the sign of $\lambda$ is determined by the conformal structure.

Considerable progress has been made on the existence problem for solutions of (\ref{PB4}); see for example \cite{ChangYangAnnals}, \cite{DM}, \cite{LLPADVMATH}, and references therein.
Our interest in this paper is dimensions $n \geq 5$, where the lack of a maximum principle (since the equation is higher order) presents an obvious difficulty when seeking positive solution of (\ref{PQ}).  Consequently, the existence theory is far less developed.  Note that for \eqref{PB4} no sign condition on
$w$ is required.

There are some results in special geometric settings.  Djadli-Hebey-Ledoux \cite{DHL} studied the optimal constant in the Sobolev embedding $W^{2,2} \hookrightarrow L^{\frac{2n}{n-4}}$
when $n \geq 5$. As a corollary of their analysis they proved some compactness results for solutions of (\ref{BPE}) assuming a size condition on $\lambda$, and that $P_g$ has constant coefficients (which holds, for example, if $(M^n,g)$ is an Einstein metric).  The assumption of constant coefficients allowed them to factor $P$ into the product of two second order operators, then apply the standard maximum principle (see also \cite{VV93}). Esposito-Robert \cite{ER} were able to find solutions to the PDE
$$
  P_g u = \lambda |u|^{\frac{8}{n-4}} u
$$
in dimension $n \geq 8$ for non-locally conformally flat manifolds, in the spirit of
\cite{Aub76}, but with no information on their sign.

 The first general existence result for (\ref{BPE}) was due to Qing-Raske \cite{QR}.  They considered locally conformally flat manifolds of positive scalar curvature, which allowed them to appeal to the work of Schoen-Yau \cite{SY} to lift the metric to a domain in the sphere via the developing map.  Assuming the Poincar\'e exponent is less than $(n-  4)/2$, they proved the existence of a
positive solution to the Paneitz-Branson equation with $\lambda  > 0$.  Hebey-Robert \cite{HR} also considered the locally conformally flat case with positive scalar curvature, and assumed in addition that the Paneitz operator and its Green function were positive.  They showed that when the Green's function satisfies a positive mass theorem, then the space of solutions to (\ref{BPE}) is compact.  Later, Humbert-Raulot \cite{HR2} verified the positive mass result (see Theorem \ref{PMT} is Section \ref{SMPSec}).  Collectively, the work of Hebey-Robert and Humbert-Raulot removed the topological assumption of Qing-Raske on the Poincar\'e exponent, but replaced it with strong positivity assumptions.

Our goal in this paper is to show that one can prove a maximum principle for $P$ and existence of solutions to (\ref{BPE}) under considerably weaker positivity assumptions. The conditions we impose are the following:  \\
\begin{align} \label{POS}
\left\{ \begin{array}{lll} \mbox{ $Q_g$ is semi-positive: $Q_g \geq 0$ and $Q_g > 0$ somewhere;}  \\
\\
\mbox{ the scalar curvature $R_g \geq 0$. }
\end{array}
\right.
\end{align}

The first main result of the paper is  \\

\noindent {\bf Theorem A.}  (See Theorem \ref{SMP} below) {\em   Let $(M^n,g)$  be a closed Riemannian manifold of dimension $n \geq 5$ satisfying (\ref{POS}). If $u \in C^4$ satisfies
\begin{align*}
P_g u \geq 0,
\end{align*}
then either $u > 0$ or $u \equiv 0$ on $M^n$.}
\vskip.2in

Theorem A is proved in Section \ref{SMPSec}, where we also show that (\ref{POS}) implies positivity of the Paneitz operator: \\

\noindent {\bf Proposition B.} (See Proposition \ref{PPosProp} below).  {\em   Let $(M^n,g)$  be a closed Riemannian manifold of dimension $n \geq 5$ satisfying (\ref{POS}).
Then the Paneitz operator is positive:
\begin{align} \label{Pposintro}
\int \phi P_{g}\phi\ dv \geq \mu(g)\int \phi^2\ dv,
\end{align}
with $\mu(g) > 0$. }
\vskip.2in

The proof is a simple extension of \cite{Gur1}, which considered the four-dimensional case.  Since $P_g > 0$, given any $p \in M^n$ the Green's function with pole at $p$, denoted $G_p$, exists.  As a corollary of Theorem A we have the positivity of $G_p$: \\

\noindent {\bf Proposition C.} (See Proposition \ref{PGF}) {\em   Let $(M^n,g)$  be a closed Riemannian manifold of dimension $n \geq 5$ satisfying (\ref{POS}).
If $G_p$ denotes the Green's function of the Paneitz operator with pole at $p \in M^n$, then $G_p > 0$
on $M^n \setminus \{p\}$.} \vskip.2in

Armed with Theorem A, we then address the question of existence of solutions to (\ref{BPE}). Given a Riemannian metric $g_0$ satisfying the positivity assumptions (\ref{POS}), we introduce a non-local flow whose stationary points are solutions of (\ref{BPE}) with $\lambda > 0$:
\begin{align} \label{uflowintro}
\left\{ \begin{array}{lll} \displaystyle \frac{\partial u}{\partial t} =  - u + \mu P_{g_0}^{-1} \big( |u|^{\frac{n+4}{n-4}} \big), \\
\\
u(\cdot, 0) = 1,
\end{array}
\right.
\end{align}
where
\begin{align} \label{mudef}
\mu = \frac{ \int u P_{g_0} u\ dv_0 }{\int |u|^{\frac{2n}{n-4}} dv_0 }.
\end{align}

Using the strong maximum principle and some elementary integral estimates, we show in Section \ref{flowSec} that the flow (\ref{uflowintro})  has a positive solution $u$ for all time $t \geq 0$.  We also show (see Section \ref{SSVar}) that the flow has a variational structure.  An important consequence of this fact is the monotonicity of the conformal volume:
\begin{align*}
\frac{d}{dt} Vol(g) = \frac{d}{dt} \int u^{\frac{2n}{n-4}}\ dv_0 \geq 0.
\end{align*}
This monotonicity property also implies the following space-time estimate:
\begin{align*}
&\int_0^{\infty} \Big( \int_{M^n} \big| - u + \mu P_{g_0}^{-1} \big( u^{\frac{n+4}{n-4}} \big)\big|^{\frac{2n}{n-4}}\ dv_0 \Big)^{\frac{n-4}{n}}\ dt < \infty.
\end{align*}
Using these facts, it is possible to choose a sequence of times $t_j \nearrow +\infty$ so that the sequence $u_j = u(\cdot, t_j)$ has a weak limit which is a solution of the
$Q$-curvature equation.

To rule out trivial limits, in Sections \ref{Sechighd}, \ref{Seclowd}, and \ref{converge} we show that it is possible to choose an initial metric in the conformal class of $g_0$ for which the solution of the flow satisfies
\begin{align*}
\int u^2\ dv_0 \geq \epsilon_0 > 0
\end{align*}
for all time.  The idea is to construct a test function whose Paneitz-Sobolev quotient is strictly less than the Euclidean value, and use this test function to define an initial conformal metric satisfying the positivity assumptions $(i) - (ii)$ above.

When the dimension is $n = 5,6,$ or $7$ or the manifold is locally conformally flat (LCF), the construction of initial data relies on a local expansion on the Green's function of the Paneitz operator.  This is proved in Section \ref{GreenF}, where we also prove a positive mass theorem.  The positive mass result extends
the version of Humbert-Raulot in \cite{HR2}, which they proved in the LCF setting (see Proposition \ref{p:green} and Theorem \ref{PMT}).  When $n \geq 8$ and the metric is not locally conformally flat we exploit instead some estimates of Esposito-Robert in \cite{ER}.  In all cases, we need to find \underline{positive} test functions with
semi-positive Q-curvature, and the strong maximum principle is crucial in this construction.

Finally, in Section \ref{converge} we show that the flow converges (up to choosing a suitable sequence of times) to a solution of the $Q$-curvature equation: \\

\noindent {\bf Theorem D.} (See Theorem \ref{flowcon}) {\em  Let $(M^n,g)$  be a closed Riemannian manifold of dimension $n \geq 5$ satisfying (\ref{POS}).
Then there is a conformal metric $h = u^{\frac{4}{n-4}}g$ with positive scalar curvature and constant positive $Q$-curvature. } \vskip.2in

\noindent {\bf Remarks.}  \\

\begin{enumerate}

\item After circulating a preliminary version of this manuscript, it was pointed out to us by E. Hebey and F. Robert that the maximum principle of Theorem A can be combined with compactness results in the literature, along with our positive mass result (Theorem \ref{PMT}), to give a proof of Theorem D by direct variational methods.  When the dimension $n \geq 8$ and $(M^n,g)$ is not locally conformally flat, one can use the expansions in Esposito-Robert \cite{ER} together with Proposition 4.1 and Theorem 5.2 in Robert's unpublished notes \cite{Rob09} to obtain existence.  When $n = 5, 6$ or $7$, we construct the necessary test functions to conclude compactness in Proposition \ref{p:testlowd}.  When $(M^n,g)$ is locally conformally flat, Theorem A and Propositions B and C imply that the Paneitz operator is ``strongly positive'' in the sense of Hebey-Robert \cite{HR}, and their result provides the necessary compactness theory (see also the comment at the end of their paper regarding the subcritical equation).  In particular, this implies that any conformal class of metrics which admits a metric with positive $Q$-curvature and positive scalar curvature also admits a minimizer of the total $Q$-curvature functional (with the same positivity conditions).  \\

\item Since this paper was submitted a number of preprints have appeared studying the $Q$-curvature in various settings; see \cite{HangYang1}, \cite{HangYang2}, \cite{HangYang3}, and \cite{CaseChang1}.  In particular, in \cite{HangYang2}, \cite{HangYang3} the authors have improved our result by weakening the assumption on the scalar curvature; positive Yamabe invariant is sufficient.   \\

\end{enumerate}

We conclude the Introduction by explaining how the flow (\ref{uflowintro}) is precisely the $W^{2,2}$-gradient flow for normalized total $Q$-curvature (up to a dimensional constant).  We remark that Baird-Fardoun-Regbaoui considered a non-local flow for the $Q$-curvature in four dimensions (see \cite{BFR}).  While their flow differs, some of their ideas inspired our approach.

Given a Riemannian manifold of dimension $n \geq 5$, if the Paneitz operator $P_g > 0$ then as $P$ is self-adjoint we can define the $W^{2,2}$ inner product by
\begin{align} \label{W22def} \begin{split}
\langle \phi , \psi \rangle_{W^{2,2}(g)} &= \int (P_g \phi) \psi \ dv_g \\
& \hskip-.5in = \int \big[ (\Delta_g \phi)(\Delta_g \psi) - 4 A_g(\nabla \phi,\nabla \psi) + (n-2) \sigma_1(A_g) g(\nabla \phi,\nabla \psi) + \frac{n-4}{2} Q_g \phi \psi \big] dv_g,
\end{split}
\end{align}
which induces the $W^{2,2}$-norm.  Denote the normalized total $Q$-curvature by
\begin{align} \label{totalQ}
\mathcal{Q}[g] = Vol(g)^{-\frac{n-4}{n}} \int Q_g\ dv_g.
\end{align}
By standard variational formulas, if $g' = \phi \, g$ is an infinitesimal conformal variation of a metric, then the variation of $\mathcal{Q}$ is given by
\begin{align} \label{Qpdef}
\mathcal{Q}'(g)\phi = \frac{n-4}{2} \int \phi \big( Q_g - \overline{Q}_g \big)\ dv_g,
\end{align}
where $\overline{Q}_g$ is the mean value of $Q$.  Since $P_g$ is invertible,
\begin{align} \label{Qp2} \begin{split}
\mathcal{Q}'(g)\phi &= \frac{n-4}{2} \int \phi P_g \big( P_g^{-1} \big( Q_g - \overline{Q}_g \big)\big)\ dv_g \\
&= \frac{n-4}{2} \int (P_g \phi) \big( P_g^{-1} \big( Q_g - \overline{Q}_g \big)\big)\ dv_g \\
&= \frac{n-4}{2} \big\langle \phi, P_g^{-1} \big( Q_g - \overline{Q}_g \big) \big\rangle_{W^{2,2}}.
\end{split}
\end{align}
Therefore, the negative $W^{2,2}$-gradient flow for the total $Q$-curvature is
\begin{align} \label{gf1}
\frac{\partial}{\partial t} g = - \frac{n-4}{2} P_g^{-1} \big( Q_g - \overline{Q}_g \big) \cdot g.
\end{align}

To see that (\ref{gf1}) is equivalent to our flow, write
\begin{align} \label{ID}
g = u^{\frac{4}{n-4}}g_0.
\end{align}
Using the conformal transformation law for the
$Q$-curvature we find
\begin{align} \label{ID1} \begin{split}
Q_g &= \frac{2}{n-4} u^{\frac{n+4}{n-4}} P_{g_0} u, \\
\overline{Q}_g &= \frac{2}{n-4} \frac{ \int u P_{g_0}u \ dv_0}{\int u^{\frac{2n}{n-4}}\ dv_0} = \frac{2}{n-4}\mu.
\end{split}
\end{align}
Also, by the conformal covariance of the Paneitz operator, its inverse is also covariant:
\begin{align} \label{ID2}
P_g^{-1} = u^{-1} P_{g_0}^{-1} \big( u^{\frac{n+4}{n-4}} \cdot \big).
\end{align}
Therefore, using (\ref{ID}), (\ref{ID1}), and (\ref{ID2}), we can rewrite (\ref{gf1}) as
\begin{align} \label{gf2}
\frac{\partial}{\partial t}u = \frac{n-4}{4} \big\{ -u + \mu  P_{g_0}^{-1}(u^{\frac{n+4}{n-4}}) \big\},
\end{align}
which only differs from our flow by the dimensional constant.

\

\noindent {\bf Acknowledgements.}  The authors would like to thank Sun-Yung Alice Chang for her careful 
reading of the original manuscript, and for suggesting revisions that improved the exposition in several places.

A.M. has been supported by the project FIRB-IDEAS {\em Analysis and Beyond},
 by the PRIN project {\em Variational Methods and Partial Differential Equations} and by the University
 of Warwick. M.J.G. is supported in part by the NSF grant DMS-1206661.

%%%%%%%%%%%%%%%%%%%%%%%%%%%%%%%%%%%%%%%%%%%%%%%%%%%%%%%%%%%%%%%%%%%%%%%%%%%%%%%%%%%%%%%%%%%
\section{The Paneitz operator and its Green's function} \label{GreenF}
%%%%%%%%%%%%%%%%%%%%%%%%%%%%%%%%%%%%%%%%%%%%%%%%%%%%%%%%%%%%%%%%%%%%%%%%%%%%%%%%%%%%%%%%%%%

In this section we prove various properties of the Paneitz operator and its Green's function that will be used throughout the paper.  \\

%%%%%%%%%%%%%%%%%%%%%%%%%%%%%
\subsection{Positivity of Paneitz operator and the Strong Maximum Principle} \label{SMPSec}  We begin with two results on the Paneitz operator:
a comparison principle, and a coercivity estimate.  We also prove a technical lemma; it shows that a metric with semi-positive $Q$-curvature and non-negative scalar curvature
 must have positive scalar curvature.  The proof is a simple application of the maximum principle, and a similar idea will be
 used elsewhere in the paper.  We first state the technical lemma:   \\

 \begin{lem} \label{RposLemma}    Let $(M^n,g)$  be a closed Riemannian manifold of dimension $n \geq 5$.  Assume  \\

\noindent $(i)$ $Q_g$ is semi-positive, i.e., $Q_g \geq 0$ and $Q_g > 0$ somewhere; \\

\noindent $(ii)$ The scalar curvature $R_g \geq 0$. \\

Then the scalar curvature is strictly positive: $R_g > 0$.
\end{lem}

\begin{pf}  By (\ref{Qdef}) the $Q$-curvature can be expressed as
\begin{align} \label{Qn}
Q_g = - \frac{1}{2(n-1)} \Delta_{g} R_g + c_1(n) R_g^2 - c_2(n) |Ric(g)|^2,
\end{align}
where $c_1(n), c_2(n) > 0$.  Since $Q_g$ is non-negative, it follows that
\begin{align} \label{MP1}
\frac{1}{2(n-1)} \Delta_{g} R_g \leq c_1(n) R_g^2.
\end{align}
By the strong maximum principle, either $R_g > 0$ or $R_g \equiv 0$.  In the latter case, by (\ref{Qn}) we would have
\begin{align} \label{Qneg}
Q_g =  - c_2(n) |Ric(g)|^2 \leq 0,
\end{align}
which is a contradiction.
\end{pf}
\vskip.2in

We now prove Theorem A of the Introduction:  \\

\begin{thm}  \label{SMP} Let $(M^n,g)$  be a closed Riemannian manifold of dimension $n \geq 5$.  Assume  \\

\noindent $(i)$ $Q_g$ is semi-positive, \\

\noindent $(ii)$ $R_g \geq 0$. \\

If $u \in C^4$ satisfies
\begin{align} \label{Pup}
P_g u \geq 0,
\end{align}
then either $u > 0$ or $u \equiv 0$ on $M^n$.

Moreover, if $u > 0$, then $h = u^{\frac{4}{n-4}}g$ is a metric with non-negative $Q$-curvature and positive scalar curvature
\end{thm}

\begin{pf}  For $\lambda \in [0,1]$ we let
\begin{align} \label{ulam}
u_{\lambda} = (1-\lambda) + \lambda u.
\end{align}
Then $u_0 \equiv 1$, while $u_1 = u$.  Assume
\begin{align} \label{uneg}
\min_{M^n} u \leq 0.
\end{align}
Define $\lambda_0 \in (0,1]$ by
\begin{align} \label{lzdef}
\lambda_0 = \min \{ \lambda \in (0,1]\ :\ \min_{M^n} u_{\lambda} = 0 \}.
\end{align}
Then for $0 < \lambda < \lambda_0$, it follows that $u_{\lambda} > 0$.  Let
\begin{align} \label{gldef}
g_{\lambda} = u_{\lambda}^{4/(n-4)}g,
\end{align}
and let $Q_{\lambda}$ denote the $Q$-curvature of $g_{\lambda}$. Note that for $0 < \lambda < \lambda_0$, we have
\begin{align} \label{Qlpos}
Q_{\lambda} \geq 0
\end{align}
and $Q_{\lambda} > 0$ somewhere. This follows from the transformation law for the $Q$-curvature:
\begin{align} \label{Qcalc} \begin{split}
Q_{\lambda} &= \frac{2}{n-4} u_{\lambda}^{-\frac{n+4}{n-4}} P_g u_{\lambda} \\
&= \frac{2}{n-4} u_{\lambda}^{-\frac{n+4}{n-4}} \Big\{ P_g\big( (1-\lambda) + \lambda u \big) \Big\} \\
&= \frac{2}{n-4} u_{\lambda}^{-\frac{n+4}{n-4}} \Big\{ (1-\lambda)P_g(1) + \lambda P_gu \Big\} \\
&= \frac{2}{n-4} u_{\lambda}^{-\frac{n+4}{n-4}} \Big\{ (1-\lambda)\frac{n-4}{2}Q_g + \lambda P_g u  \Big\} \\
&\geq (1-\lambda) Q_g u_{\lambda}^{-\frac{n+4}{n-4}}.
\end{split}
\end{align}
Since $\lambda < \lambda_0 \leq 1$ and $Q_g$ is semi-positive, it follows that $Q_{\lambda}$ is semi-positive.

Let $R_{\lambda}$ denote the scalar curvature of $g_{\lambda}$.  We also claim that for $0 \leq \lambda < \lambda_0$,
\begin{align} \label{Rlpos}
R_{\lambda} > 0.
\end{align}
This certainly holds for $\lambda = 0$; but if there is a $\lambda_1 \in (0,\lambda_0)$ with $\min R_{\lambda_1} = 0$, then this would
contradict Lemma \ref{RposLemma}.

By the formula for the transformation of the scalar curvature under a conformal change of metric,
\begin{align} \label{RU}
R_{\lambda} = u_{\lambda}^{-\frac{n}{n-4}} \Big\{ - \frac{4(n-1)}{(n-4)} \Delta_g u_{\lambda} - \frac{8(n-1)}{(n-4)^2} \frac{|\nabla_g u_{\lambda}|^2}{u_{\lambda}} + R_g u_{\lambda} \Big\}.
\end{align}
Since $R_{\lambda} > 0$, this implies $u_{\lambda}$ satisfies the differential inequality
\begin{align} \label{Upde}
\Delta_g u_{\lambda} \leq \frac{(n-4)}{4(n-1)} R_g u_{\lambda}.
\end{align}
Taking the limit as $\lambda \nearrow \lambda_0$, this also holds for $\lambda = \lambda_0$.
By the strong maximum principle, (\ref{lzdef}) and (\ref{Upde}) imply $u_{\lambda_0} \equiv 0$.  If $\lambda_0 = 1$, then we are done.
Therefore, assume $\lambda_0 \in (0,1)$.  It follows from (\ref{ulam}) that
\begin{align*}
u = - \frac{(1- \lambda_0)}{\lambda_0},
\end{align*}
hence
\begin{align*}
P_g u =  - \big(\frac{n-4}{2}\big) \frac{(1- \lambda_0)}{\lambda_0} Q_g.
\end{align*}
Since by assumption $Q_g > 0$ somewhere, this contradicts $P_g u \geq 0$.  We conclude that $u \equiv 0$ or $u > 0$.

If $u > 0$, then the metric $h = u^{\frac{4}{n-4}}g$ is well defined and has non-negative $Q$-curvature. Once again, we can define the family of functions $\{ u_{\lambda} \}$ as
in (\ref{ulam}) and the metrics $g_{\lambda}$ as in (\ref{gldef}).  Then the scalar curvature of $g_{\lambda}$ satisfies (\ref{RU}), and by the strong maximum principle it follows
that either $R_{\lambda} > 0$ or $R_{\lambda} \equiv 0$.  Recall by Lemma \ref{RposLemma} that $R_g > 0$.  Therefore, we cannot have $R_{\lambda} \equiv 0$, since a conformal class which admits a metric of positive scalar curvature cannot admit a scalar-flat metric.  It follows that $R_{\lambda} > 0$ for all $\lambda \in [0,1]$.
\end{pf}

We now show that the positivity assumptions of the preceding theorem imply the positivity of the Paneitz operator. This is easy to prove in dimensions $n \geq 6$, but
 for $n = 5$ we need to adapt the idea of the $n=4$ case appearing in \cite{Gur1}. \\

\begin{pro}  \label{PPosProp} Under the assumptions of Theorem \ref{SMP} the Paneitz operator is positive:
there exists $\mu(g) > 0$ such that
\begin{align} \label{Ppos}
\int \phi P_{g}\phi\ dv \geq \mu(g) \int \phi^2 dv.
\end{align}
Consequently, the Paneitz-Sobolev constant is also positive:
\begin{align} \label{qdef}
q_0(M^n,g) \equiv \inf_{\phi \in W^{2,2} \setminus \{0\}} \frac{ \displaystyle \int \phi P_{g} \phi\ dv }{ \displaystyle \big( \int |\phi|^{\frac{2n}{n-4}}\ dv \big)^{\frac{n-4}{n}} } > 0.
\end{align}
\end{pro}
\vskip.1in

\begin{pf}  By (\ref{Pdef}),
\begin{align} \label{DF1}
\int \phi P \phi\ dv = \int \big\{ (\Delta \phi)^2 - 4 A(\nabla \phi,\nabla \phi) + (n-2) \sigma_1(A)|\nabla \phi|^2 + \frac{n-4}{2}Q \phi^2 \big\}\ dv,
\end{align}
where we have omitted the subscript $g$. There are two cases to consider: $n = 5$, and $n \geq 6$.  In the latter case we use the integrated Bochner formula
\begin{align*}
\int (\Delta \phi)^2\ dv &= \int |\nabla^2 \phi|^2\ dv + \int Ric(\nabla \phi,\nabla \phi)\ dv \\
&= \int |\nabla^2 \phi|^2\ dv + (n-2) \int A(\nabla \phi,\nabla \phi)\ dv + \int \sigma_1(A)|\nabla \phi|^2\ dv,
\end{align*}
which gives
\begin{align} \label{Boch}
 \int - 4 A(\nabla \phi,\nabla \phi)\ dv = \int \big\{ -\frac{4}{n-2}(\Delta \phi)^2 + \frac{4}{n-2}|\nabla \phi|^2 + \frac{4}{n-2}\sigma_1(A)|\nabla \phi|^2 \big\}\ dv.
 \end{align}
 Substituting this into (\ref{DF1}) we find
 \begin{align} \label{DF2}
 \int \phi P \phi\ dv = \int \Big\{ \frac{n-6}{n-2} (\Delta \phi)^2 + \frac{4}{n-2}|\nabla^2 \phi|^2  + \frac{(n-2)^2 + 4}{n-2} \sigma_1(A)|\nabla \phi|^2 + \frac{n-4}{2}Q \phi^2 \Big\}\ dv.
 \end{align}
Consequently, when the dimension $n \geq 6$ the positivity of $P$ follows.

When $n = 5$ we need to adapt the argument for the four-dimensional case in \cite{Gur1}.  First, when $n=5$ we note that
\begin{align}  \label{DF5}
\int \phi P\phi\ dv = \int (\Delta \phi)^2\ dv - 4 \int A(\nabla \phi,\nabla \phi)\ dv + 3 \int \sigma_1(A)|\nabla \phi|^2\ dv + \frac{1}{2}\int Q \phi^2\ dv,
\end{align}
while the $Q$-curvature is given by
\begin{align}  \label{Q5}
0 \leq Q = -\Delta \sigma_1(A) - 2 |A|^2 + \frac{5}{2} \sigma_1(A)^2.
\end{align}

Consider the second term on the right-hand side of (\ref{DF5}).  Since by Lemma \ref{RposLemma} the scalar curvature is positive, using the arithmetic/geometric mean inequality (AGM) we estimate
\begin{align*}
4 A(\nabla \phi, \nabla \phi) \leq 2 \frac{|A|^2}{\sigma_1(A)}|\nabla \phi|^2 + 2 \sigma_1(A) |\nabla \phi|^2.
\end{align*}
By (\ref{Q5}),
\begin{align}  \label{leq5}
2 \frac{|A|^2}{\sigma_1(A)}|\nabla \phi|^2 \leq - \frac{\Delta \sigma_1(A)}{\sigma_1(A)} |\nabla \phi|^2 + \frac{5}{2}\sigma_1(A)|\nabla \phi|^2,
\end{align}
hence
\begin{align} \label{leq51}
4 \int A(\nabla \phi, \nabla \phi) \leq - \int \frac{\Delta \sigma_1(A)}{\sigma_1(A)} |\nabla \phi|^2\ dv  + \frac{9}{2} \int \sigma_1(A)|\nabla \phi|^2\ dv.
\end{align}
For the first term on the right, we integrate by parts and use the AGM inequality to get
\begin{align} \label{leq52} \begin{split}
- \int \frac{\Delta \sigma_1(A)}{\sigma_1(A)} |\nabla \phi|^2\ dv &= \int \Big\{  - \frac{|\nabla \sigma_1(A)|^2}{\sigma_1(A)^2}|\nabla \phi|^2 + \big\langle \frac{\nabla \sigma_1(A)}{\sigma_1(A)}, \nabla |\nabla \phi|^2 \big\rangle\Big\} dv \\
&= \int \Big\{  - \frac{|\nabla \sigma_1(A)|^2}{\sigma_1(A)^2}|\nabla \phi|^2 + 2 \nabla^2 \phi \big( \frac{\nabla \sigma_1(A)}{\sigma_1(A)}, \nabla \phi \big) \Big\}dv \\
&\leq \int \Big\{  - \frac{|\nabla \sigma_1(A)|^2}{\sigma_1(A)^2}|\nabla \phi|^2 + 2 |\nabla^2 \phi|\frac{|\nabla \sigma_1(A)|}{\sigma_1(A)}|\nabla \phi| \Big\}dv  \\
&\leq \int \Big\{  - \frac{|\nabla \sigma_1(A)|^2}{\sigma_1(A)^2}|\nabla \phi|^2 +  \frac{|\nabla \sigma_1(A)|^2}{\sigma_1(A)^2}|\nabla \phi|^2 + |\nabla^2 \phi|^2  \Big\}dv \\
&= \int |\nabla^2 \phi|^2\ dv.
\end{split}
\end{align}
Substituting this back into (\ref{leq51}) gives
\begin{align} \label{leq53}
4 \int A(\nabla \phi, \nabla \phi) \leq  \int |\nabla^2 \phi|^2\ dv  + \frac{9}{2} \int \sigma_1(A)|\nabla \phi|^2\ dv.
\end{align}
In dimension five the Bochner formula gives
\begin{align*}
\int |\nabla^2 \phi|^2\ dv = \int (\Delta \phi)^2\ dv - 3 \int A(\nabla \phi,\nabla \phi)\ dv - \int \sigma_1(A)|\nabla \phi|^2\ dv,
\end{align*}
and substituting this into (\ref{leq53}) we arrive at
\begin{align} \label{leq54}
4 \int A(\nabla \phi, \nabla \phi) \leq \int (\Delta \phi)^2\ dv - 3 \int A(\nabla \phi,\nabla \phi)\ dv + \frac{7}{2} \int \sigma_1(A)|\nabla \phi|^2\ dv.
\end{align}
Combining the Schouten tensor terms we have
\begin{align} \label{leq55}
7 \int A(\nabla \phi, \nabla \phi) \leq \int (\Delta \phi)^2\ dv + \frac{7}{2} \int \sigma_1(A)|\nabla \phi|^2\ dv,
\end{align}
hence
\begin{align} \label{leq56}
4 \int A(\nabla \phi, \nabla \phi) \leq \frac{4}{7} \int (\Delta \phi)^2\ dv + 2 \int \sigma_1(A)|\nabla \phi|^2\ dv,
\end{align}
or
\begin{align} \label{leq57}
- 4 \int A(\nabla \phi, \nabla \phi) \geq - \frac{4}{7} \int (\Delta \phi)^2\ dv - 2 \int \sigma_1(A)|\nabla \phi|^2\ dv.
\end{align}
Finally, substituting this into (\ref{DF5}) gives
\begin{align}  \label{DF55} \begin{split}
\int \phi P\phi\ dv &= \int (\Delta \phi)^2\ dv - 4 \int A(\nabla \phi,\nabla \phi)\ dv + 3 \int \sigma_1(A)|\nabla \phi|^2\ dv + \frac{1}{2}\int Q \phi^2\ dv \\
&\geq \frac{3}{7} \int (\Delta \phi)^2\ dv  + \int \sigma_1(A)|\nabla \phi|^2\ dv + \frac{1}{2}\int Q \phi^2\ dv,
\end{split}
\end{align}
and the positivity of $P$ follows.
\end{pf}
\vskip.2in

From Proposition \ref{PPosProp} we conclude that under the assumptions of Lemma \ref{RposLemma}, for any $p \in M^n$ the Green's function of the Paneitz operator $G_p$ exists, satisfying
\begin{align} \label{Gnorm}
P_g G_p = \delta_p,
\end{align}
where $\delta_p$ is the Dirac mass at $p$.  We now prove Proposition C of the Introduction: \vskip.1in

%We now come to one of the main results of this section: \vskip.2in

\begin{pro} \label{PGF} Suppose $(M^n,g)$ satisfies the same assumptions of Theorem \ref{SMP}.
If $G_p$ denotes the Green's function of the Paneitz operator with pole at $p \in M^n$, then $G_p > 0$
on $M^n \setminus \{p\}$.
\end{pro}

\begin{pf}  Consider a sequence of continuous functions $f_j$ on $M$ which are non-negative, whose
supports shrink to $\{p\}$, and such that
$$
  \int_M f_j dv = 1 \qquad \quad \hbox{ for all } j.
$$
Then $f_j \rightharpoonup \d_p$ in the sense of distributions. If $G_j$ is the solution to
$$
  P_g G_j = f_j,
$$
it is easy to show that
$$
  G_j \to G_p \qquad \hbox{ in } C^4_{loc}(M^n \setminus \{p\}).
$$
By Theorem \ref{SMP} one has $G_j > 0$ on $M^n$, which immediately implies that
$$
  G_p \geq 0 \qquad \quad \hbox{ on } M^n \setminus \{p\}.
$$

Suppose there exists $x_0 \neq p$ such that $G_p(x_0) = 0$, and consider the sequence of conformal metrics $g_j = G_j^{\frac{4}{n-4}} g$.  By construction $P_g G_j \geq 0$,
hence by Theorem \ref{SMP} the metrics $g_j$ have positive scalar curvature and semi-positive $Q$-curvature.  It follows that the scalar curvature of $g_j$ satisfies
\begin{align} \label{MPj}
\frac{1}{2(n-1)} \Delta_{g_j} R_{g_j} \leq c_1(n) R_{g_j}^2.
\end{align}
Also, arguing as we did in the proof of Lemma \ref{RposLemma} (see \eqref{Upde}), $G_j$ satisfies the
 differential inequality
 \begin{align} \label{Updej}
 \Delta_g G_j \leq \frac{(n-4)}{4(n-1)} R_g G_j \qquad \quad \hbox{ on } M^n.
 \end{align}
Passing to the limit $j \to \infty$ on $M^n \setminus \{p\}$ we have
$$
  \Delta_g G_p \leq \frac{(n-4)}{4(n-1)} R_g G_p.
$$
By the strong maximum principle, $G_p(x_0) = 0$ implies $G_p \equiv 0$, which is a contradiction.
\end{pf}

%%%%%%%%%

\subsection{Regularity of the Green's function}  \label{GreenApp}  Our next results concern the behavior of the Green's function near the pole. We will show that if the
 dimension is $5, 6$ or $7$, or if the manifold is locally conformally flat, then in conformal normal
coordinates the Green's function of the Paneitz operator is equal to sum of the fundamental
solution of the bi-laplace equation and a weighted Lipschitz function:   \\

\begin{pro}\label{p:green}  Let $(M^n,g)$ be a closed Riemannian manifold, satisfying the assumptions of Lemma \ref{RposLemma}:  \\

\noindent $(i)$ $Q_g$ is semi-positive, and \\

\noindent $(ii)$ $R_g \geq 0$. \\

In addition, assume one of the following holds: \\

\noindent $\bullet$ The dimension $n = 5, 6, $ or $7$; or \\

\noindent $\bullet$ $(M^n,g)$ is locally conformally flat and $n \geq 5$.  \\

For $p \in M$, consider the conformal normal coordinates centered at $p$ constructed in \cite{LP}
with conformal metric $\tilde{g}$. Then, if
$G_p(x)$ is the Green's function for the Paneitz operator with pole at $p$, there exists
a constant $\alpha$ such that in conformal normal coordinates,
\begin{align} \label{Gexp}
  G_p(x) = \frac{c_n}{d_{\tilde{g}}(x,p)^{n-4}} + \alpha + O^{(4)}(r),
\end{align}
where $c_n = \frac{1}{(n-2)(n-4)\omega_{n-1}}$, $\o_{n-1} = |S^{n-1}|$, and $f = O^{(k)}(r^m)$ denotes any quantity satisfying
\begin{align*}
|\nabla^j f(x)| \leq C_j r^{m - j}
\end{align*}
for $1 \leq j \leq k$, where $r = |x| = d_{\tilde{g}}(x,p)$.
\end{pro}

\begin{pf}  In the locally conformally flat case, one can conformally change and use Euclidean coordinates near $p$, and the
expansion (\ref{Gexp}) appears in \cite{HR2}.  For the non-LCF cases we will use the classical method of the parametrix; namely we start
with functions which properly approximate $G_p$ and then use elliptic regularity theory.  We begin with some preliminary lemmas. \\

\begin{lem}\label{l:d1d2d4u}
In conformal normal coordinates, if $u$ is a radial function then one has the following expansions:
\begin{align} \label{HessUx}
  \n_i \n_j u = \frac{x_i x_j}{r^2} u'' - \frac{x_i x_j}{r^3} u' + \frac{\d_{ij}}{r} u'
  + O(r) |u'|;
\end{align}
\begin{align} \label{LapUx}
  \D_{\tilde{g}} u = u'' + \frac{n-1}{r} u' + O''(r^{N-1}) u';
\end{align}
\begin{align} \label{BiLapx}
   \D^2_{\tilde{g}} u = \D ^2_{0} u + O(r^{N-1}) u''' + O(r^{N-2}) u'' + O(r^{N-3}) u',
\end{align}
where $N \geq 5$ and $\Delta_0$ denotes the Euclidean laplacian.
\end{lem}

\begin{pf}
Let $\{ x^i \}$ denote conformal normal coordinates associated with the metric $\tilde{g}$, and let $\{ r, \vartheta^{\alpha} \}$ denote the corresponding
polar coordinates, where $r = |x|$ and $\{ \vartheta^{\alpha}\}$ are coordinates on the unit sphere.   We let $\tilde{g} = \tilde{g}_{ij}$ denote the matrix of components of $\tilde{g}$ with respect to the $\{ x^i \}$ coordinates,
and $\tilde{g}^{\prime} = \tilde{g}_{\alpha \beta}^{\prime}$ the components of $\tilde{g}$ with respect to the polar coordinate system.  It follows
that
\begin{align} \label{gtogp}
 \sqrt{\det \tilde{g}^{\prime}} = r^{n-1} \sqrt{\det \tilde{g}}.
 \end{align}
If $u$ is radial, then
\begin{align} \label{DGu} \begin{split}
  \D u(r) &= \frac{1}{\sqrt{\det \tilde{g}^{\prime}}} \pa_r \left( \sqrt{\det {\tilde{g}^{\prime}}} \, \pa _r u \right) \\
  &= u'' + \partial_r \big(\log \sqrt{\det \tilde{g}^{\prime}}\big) u' \\
  &= u'' + \partial_r \big(\log r^{n-1} \sqrt{\det \tilde{g}} \big) u' \\
  &= u'' + \frac{n-1}{r} u' + u' \partial_r \log \sqrt{\det \tilde{g}} \\
  &= \Delta_0 u + u' \partial_r \log \sqrt{\det \tilde{g}}. \\
  \end{split}
  \end{align}
In conformal normal coordinates (see Theorem 5.1 of \cite{LP}) the determinant of $\tilde{g}$ approaches
$1$ smoothly at the origin at order $N$, where $N \geq 5$, and in particular one has
\begin{align} \label{gN}
\det \tilde{g} = 1 + O^{(3)}(r^N).
\end{align}
Therefore,
\begin{align} \label{DlogD}
\partial_r \log \sqrt{\det \tilde{g}} = O''(r^{N-1}).
\end{align}
Substituting into (\ref{DGu}), we arrive at (\ref{LapUx}). The formula (\ref{BiLapx}) for the bi-laplacian follows immediately.

Recall  that, in normal coordinates
\begin{align} \label{gnormexp}
  {\tilde{g}}_{ij} = \d_{ij} - \frac{1}{3} R_{i \a j \b} x^\a x^\b + O^{(4)}(r^3),
\end{align}
where $R_{*}$ denotes the curvature tensor (with respect to $\tilde{g}$) evaluated at $p$.
%
%
%Using the decomposition of the curvature tensor,
%\begin{eqnarray*}
%W_{ijkl} & = & R_{ijkl} - \frac{1}{n-2} \left( R_{ik}\tilde{g}_{jl} -
%R_{il}\tilde{g}_{jk} + R_{jl}\tilde{g}_{ik} - R_{jk}\tilde{g}_{il} \right)
% \\ & + & \frac{R}{(n-1)(n-2)}(\tilde{g}_{jl}\tilde{g}_{ik} - \tilde{g}_{jk}\tilde{g}_{il}),
%\end{eqnarray*}
%and the fact that the Ricci tensor and the scalar curvature vanish at $p$, we can rewrite (\ref{gnormexp}) as
%$$
%    \tilde{g}_{ij} = \d_{ij} + \frac{1}{3} W_{i \a j \b} x^\a x^\b + O^{(4)}(r^3).
%$$
%Therefore,
%\begin{eqnarray*}
%   \partial_i \tilde{g}_{jm} & = & \partial_i \left[\d_{jm} - \frac{1}{3} W_{i \a  m \b} x^\a x^b +
%   O(r^3) \right] \\ & = & - \frac{1}{3} W_{i \a m \b} \d_i^\a x^\b - \frac{1}{3}
%   W_{i \a m \b} x^\a \d_i^\b + O'''(r^2) = O'''(r^2).
%\end{eqnarray*}
As
$$
  \G^k_{ij} = \frac{1}{2} g^{km} \left[ \partial_i \tilde{g}_{jm} + \partial_j \tilde{g}_{im}
  - \partial_m \tilde{g}_{ij} \right],
$$
we deduce that
$$
   | \G^k_{ij}(x) | = O'''(r); \qquad \qquad |\partial_l \G^k_{ij}(x)| = O''(1).
$$
This implies that
$$
  \n_i \n_j u = \partial_i \partial_j u   + O(r) |u'|.
$$
As $u$ is radial, we obtain the conclusion.
\end{pf}

\noindent {\bf Remark.}  In the estimates that follow we will only need the order of flatness in (\ref{gN}) to be $N = 5$.  Therefore, we will assume from now on that $N \geq 5$ is fixed.
\\

\begin{lem}\label{l:expgeom}
In conformal normal coordinates one has the following expansions for the
Schouten tensor $A_{ij} = \frac{1}{n-2} \left[ R_{ij} - \frac{1}{2(n-1) }R_{\tilde{g}} \, \tilde{g}_{ij} \right]$ and for the $Q$-curvature:
$$
   A_{ij}(0) = 0; \qquad \qquad
   \left( \n_k A_{ij} + \n_i A_{jk} + \n_j A_{ik} \right)(0)  = 0;
$$
$$
\n_k \n _l A_{ij}(0) x^k x^l x^i x^j = - \frac{r^2}{(n-2)} \n_k \n _l \s_1(0) x^k x^l;
$$
$$
  Q = - \frac{1}{2(n-1)} \left[ - \frac{1}{6} |W|^2(0) + O(r) \right].
$$
\end{lem}

\begin{pf}
Recall that in conformal normal coordinates one has
$$
  R_{ij}(0) = 0; \qquad \qquad
     \left( \n_k R_{ij} + \n_i R_{jk} + \n_j R_{ik} \right)(0)  = 0;
$$
$$
  \left( \n_k \n _l R_{ij} + \n_l \n _i R_{jk} + \n_i \n _j R_{kl} + \n_j \n _k R_{li}
  \right)(0) = 0;
$$
$$
   R(0) = 0; \qquad \n_{\tilde{g}} R(0) = 0; \qquad \D_{\tilde{g}} R(0) = \frac{1}{6} |W|^2(0).
$$
Then the conclusion follows immediately from the definition of $A_{ij}$
and $Q$.
\end{pf}

\

\begin{lem}\label{l:panrad}
If $u$ is a radial function, then in conformal normal coordinates and conformal metric
$\tilde{g}$ one has that
\begin{eqnarray}\label{eq:exprad}\nonumber
    P_{\tilde{g}} u & = & \D ^2_{0} u + \n_k \n_l \s_1(0) x^k x^l  \mathfrak{Q}(u) +
    \frac{n-4}{24(n-1)} |W|^2(0) \, u +    O(r^3) |u''|
  \\ & + &  O(r^2) |u'| + O(r) u  + O(r^{N-1}) u''' + O(r^{N-2}) u'' + O(r^{N-3}) u',
\end{eqnarray}
where
$$
  \mathfrak{Q}(u) =  \frac{u'}{r} \left(
      \frac{2(n-1)}{(n-2)} - \frac{(n-1)(n-2)}{2} + 6 - n \right)
      - u'' \left( \frac{(n-2)}{2} + \frac{2}{(n-2)} \right).
$$
\end{lem}

\begin{pf}
Recall that
$$
  P_g u = \Delta^2_{\tilde{g}} u + \mbox{div}_{\tilde{g}} \big\{ \big( 4 A_{\tilde{g}}  - (n-2) \sigma_1(A_{\tilde{g}} )
  \tilde{g} \big)(\nabla u, \cdot) \big\} + \frac{n-4}{2} Q_{\tilde{g}} u.
$$
We consider  the term
$$
   \mbox{div}_{\tilde{g}} \big\{ \big( 4 A_{\tilde{g}}  - (n-2) \sigma_1(A_{\tilde{g}} )
     \tilde{g} \big)(\nabla u, \cdot) \big\}  = 4 A_{ij} \n_i \n_j u - (n-2) \s_1 \D_{\tilde{g}}   u + (6-n)
   \langle \n \s_1,  \n u \rangle.
$$

Using Lemma \ref{l:d1d2d4u} for the Hessian of $u$ and Lemma \ref{l:expgeom} for the
vanishing of $A_{ij}(0)$ we find that
\begin{eqnarray} \nonumber
A_{ij} \n_i \n_j u  & = & (A_{ij}(0) + \n_k A_{ij}(0) x^k + \frac{1}{2} \n_k \n_l A_{ij}(0) x^k x^l + O(r^3))
  \left( \partial^2_{ij} u + O(r) |u'| \right)  \\ & = &
  (\n_k A_{ij}(0) x^k + \frac{1}{2} \n_k \n_l A_{ij}(0) x^k x^l + O(r^3))   \\ & \times & \nonumber
  \left[ \left( \frac{\d^{ij}}{r^2} - \frac{x^i x^j}{r^3}  \right) u' + \frac{x^i x^j}{r^2} u'' + O(r) |u'|  \right]
  = I + II + III + IV + V,
\end{eqnarray}
where
$$
  I = \n_k A_{ij}(0) x^k \frac{\d^{ij}}{r^2} u'; \qquad \qquad II = \n_k A_{ij}(0)
  x^k \left( \frac{x^i x^j}{r^2} u''  - \frac{x^i x^j}{r^3}  u' \right);
$$
$$
  III =  \frac{1}{2} \n_k \n_l A_{ij}(0) x^k x^l \frac{\d^{ij}}{r^2} u'; \qquad \qquad IV =
  \frac{1}{2} \n_k \n_l A_{ij}(0) x^k x^l \left( \frac{x^i x^j}{r^2} u''  - \frac{x^i x^j}{r^3}  u' \right);
$$
$$
  V = O(r^3) \left( \partial^2_{ij} u + O(r) |u'| \right)  +
   (\n_k A_{ij}(0) x^k + \frac{1}{2} \n_k \n_l A_{ij}(0) x^k x^l) \times  O(r) |u'|.
$$
As the scalar curvature vanishes to first order at $p$ we find immediately that $I = 0$.  Also, since $II$ stays unchanged
after permutation of the indices $i,j,k$, by the second statement of Lemma \ref{l:expgeom} we find
that also $II = 0$. Turning to $III$, we have that
$$
  III = \frac{1}{2} \n_k \n_l \s_1 x^k x^l \frac{u'}{r}.
$$
Concerning  $IV$ instead, using the third identity in Lemma \ref{l:expgeom}  we find that
$$
  IV = - \frac{1}{2(n-2)} \n_k \n_l \s_1 x^k x^l \left( u'' - \frac{u'}{r} \right).
$$
Expanding then also $V$ one finds
$$
 4 A_{ij} \n_i \n_j u  = \frac{2}{n-2} \n_k \n_l \s_1(0) x^k x^l \left[ (n-1) \frac{u'}{r} - u'' \right]
 + O(r^3) |u''| + O(r^2) |u'|.
$$
Similarly, using the second assertion of Lemma \ref{l:d1d2d4u} and a Taylor expansion of the
scalar curvature one finds
$$
- (n-2)  \s_1 \D u = - \frac{n-2}{2} \n_k \n_l \s_1(0) x^k x^l
\left[ (n-1) \frac{u'}{r} + u'' \right] +  O(r^3) |u''| + O(r^2) |u'|.
$$
Furthermore
$$
  (6-n)
     \langle \n \s_1,  \n u \rangle = (6-n) \n_k \n_l \s_1(0) x^k x^l  \frac{u'}{r}
     + O(r^2) |u'|.
$$
By the third assertion of Lemma \ref{l:d1d2d4u}  and
summing all the above terms in $P_{\tilde{g}} u$ (taking into account of the expression of
$Q_{\tilde{g}}$ in Lemma \ref{l:expgeom}) one gets the conclusion.
\end{pf}

\

Using the preceding technical lemmas, we can now compute $P_{\tilde{g}} (r^{4-n})$.  By Lemma \ref{l:panrad}, one has that
\begin{align} \label{Prn}
  P_{\tilde{g}} (r^{4-n}) = \mathfrak{A}_n \d_p + \n_k \n_l \s_1(0) x^k x^l  \mathfrak{Q}_n r^{2-n} +
   \frac{n-4}{24(n-1)} |W|^2(0) \, r^{4-n}  +  O(r^{5-n}),
\end{align}
where $\mathfrak{A}_n = 2 (n-2) (n-4) |S^{n-1}|$, and where
$$
  \mathfrak{Q}_n =  (4-n) \left[ \left(
      \frac{2(n-1)}{(n-2)} - \frac{(n-1)(n-2)}{2} + 6 - n \right)
      - (3-n) \frac{(n-2)^2+4}{2(n-2)}  \right].
$$
It follows from (\ref{Prn}) that
\begin{align} \label{Errest}
  P_{\tilde{g}} \left( G_p - \frac{1}{\mathfrak{A}_n}  r^{4-n} \right) = O(r^{4-n}).
\end{align}
By elliptic regularity, if we can show that the right-hand side of (\ref{Errest}) is in $L^p$ for some $p > n/3$, then we would conclude
\begin{align*}
 G_p - \frac{1}{\mathfrak{A}_n}  r^{4-n}  \in W^{4,p} \hookrightarrow C^{1,\alpha}
\end{align*}
with $\alpha > 0$, and (\ref{Gexp}) would follow.  However, $r^{4-n} \in L^p  \ \mbox{ for } \ p < \frac{n}{n-4}$,
hence we need $p$ to satisfy
\begin{align*}
\frac{n}{3} < p < \frac{n}{n-4}.
\end{align*}
This can only hold if $n=5$ or $n = 6$; when $n = 7$ we have equality, so this is the borderline case.

When $n=7$ we can add a further correction term to study the asymptotics of $G_p$. We begin by writing the trailing terms in (\ref{Prn}) as
$$
  \n_k \n_l \s_1(0) x^k x^l  \mathfrak{Q}_n |x|^{2-n} +
     \frac{n-4}{24(n-1)} |W|^2(0) \, |x|^{4-n} = \mathfrak{B}_0 |x|^{-3}
     + \mathfrak{B}_2(\theta) |x|^{-3},
$$
where $\mathfrak{B}_0$ is a  constant and where $\mathfrak{B}_2(\theta)$
is a second spherical harmonic function (with zero average) on $S^6$, with
$\theta$ denoting the spherical coordinates.  As the second eigenvalue of the Laplace-Beltrami operator on $S^6$ is equal to 14,
using polar coordinates one can easily check that
$$
  \D^2_{0} |x| = - \frac{24}{|x|^3};   \qquad \qquad
  \D^2_{0} \left( \mathfrak{B}_2(\theta) |x| \right) = \frac{172}{|x|^3}.
$$
Therefore in conformal normal coordinates one finds that
$$
  \D \left( - \frac{1}{24} \mathfrak{B}_0 |x| + \frac{\mathfrak{B}_2(\theta)}{172} |x| \right)
  = \mathfrak{B}_0 |x|^{-3}  + \mathfrak{B}_2(\theta) |x|^{-3} + O(r^{-2}),
$$
which implies that
$$
  P_{\tilde{g}} \left(  G_p - \frac{1}{\mathfrak{A}_n}  |x|^{4-n}   + \frac{1}{24} \mathfrak{B}_0 |x| - \frac{\mathfrak{B}_2(\theta)}{172} |x| \right) = O(r^{-2}).
$$
By elliptic regularity theory and by Morrey's embedding theorems we then deduce that
the function
$$
  G_p - \frac{1}{\mathfrak{A}_n}  |x|^{4-n}   + \frac{1}{24} \mathfrak{B}_0 |x| - \frac{\mathfrak{B}_2(\theta)}{172} |x|
$$
possesses H\"older continuous derivatives which, taking Schauder's estimates into account,
implies the conclusion when $n=7$.
\end{pf}
\vskip.2in

%%%%%%%%%%%%%%%%%%%%%%%%%%%%%%%%%%%%%%%%%%%%%%%%%
\subsection{A positive mass theorem for the Paneitz operator}  We conclude this section by proving an inequality for the constant $\alpha$ in the expansion for
the Green's function in Proposition \ref{p:green}.  In the locally conformally flat case, this was proved by Humbert-Raulot in \cite{HR2}.  In fact, their
proof is easily adapted to the non-LCF case when the dimension is $5,6, $ or $7$.  \\

\begin{thm} \label{PMT}  Under the assumptions of Proposition \ref{p:green}, the constant $\alpha$ in the expansion (\ref{Gexp}) satisfies $\alpha \geq 0$, with
equality if and only if $(M^n,g)$ is conformally equivalent to the round sphere.
\end{thm}

\begin{pf}  Let $\Gamma_p$ denote the Green's function for the conformal laplacian $L = - \Delta + \frac{(n-2)}{4(n-1)}R$ with pole at $p$. As in \cite{HR2}, we consider the conformal blow-up of $g$ defined by
\begin{align} \label{Gamg}
\hat{g} = \Gamma_p^{\frac{4}{n-2}} g.
\end{align}
This defines an asymptotically flat, scalar-flat metric on $X^n = M^n \setminus \{ p \}$.  Let
\begin{align} \label{Phidef}
\Phi = \Gamma_p^{-\frac{n-4}{n-2}}G_p.
\end{align}
By the conformal covariance of the Paneitz operator, on $X^n$ we have
\begin{align*}
P_{\hat{g}} \Phi &= P_{ \Gamma_p^{\frac{4}{n-2}} g } \big( \Gamma_p^{-\frac{n-4}{n-2}}G_p \big) \\
&= \Gamma_p^{-\frac{n+4}{n-2}} P_{g} (G_p) \\
&= 0.
\end{align*}
Also, since $\hat{g}$ is scalar flat, its $Q$-curvature is given by
\begin{align} \label{Qhg}
Q_{\hat{g}} = - 2 |A(\hat{g})|^2,
\end{align}
where $A$ is the Schouten tensor.  By the formula for the Paneitz operator (\ref{Pdef}),
\begin{align} \label{PPhi}
0 = P_{\hat{g}} \Phi = \Delta_{\hat{g}}^2 \Phi + \mbox{div}_{\hat{g}} \big\{ \big( 4 A_{\hat{g}} (\nabla \Phi, \cdot) \big\} - (n-4) |A(\hat{g})|^2 \Phi.
\end{align}

Fix $\delta > 0$ small and let $B_{\delta}$ once again denote the geodesic ball centered at $p$ of radius $\delta > 0$ (as measured in the metric $g$ -- not $\hat{g}$).
As in \cite{HR2}, we integrate (\ref{PPhi}) over $M^n \setminus B_{\delta}$ and apply the divergence theorem:
\begin{align} \label{Is1} \begin{split}
0 &= \int_{M^n \setminus B_{\delta}} P_{\hat{g}} \Phi\ dv_{\hat{g}}  \\
&= \int_{M^n \setminus B_{\delta}} \Big\{ \Delta_{\hat{g}}^2 \Phi + \mbox{div}_{\hat{g}} \big\{ \big( 4 A_{\hat{g}} (\nabla \Phi, \cdot) \big\} - (n-4) |A(\hat{g})|^2 \Phi \Big\}\ dv_{\hat{g}} \\
&= \oint_{\partial B_{\delta}} \big\{ \frac{\partial}{\partial \nu} \big(\Delta_{\hat{g}} \Phi \big) + 4 A_{\hat{g}}(\nabla \Phi, \nu) \big\} dS_{\hat{g}} - (n-4) \int_{M^n \setminus B_{\delta}} |A(\hat{g})|^2 \Phi\ dv_{\hat{g}},
\end{split}
\end{align}
where $\nu$ is the (outward) normal to $\partial B_{\delta}$ in the metric $\hat{g}$.

Considering the boundary integrals, we first note that since $\hat{g}$ is scalar-flat,
\begin{align}  \label{PL}
 \frac{\partial}{\partial \nu} \big(\Delta_{\hat{g}} \Phi \big) =  - \frac{\partial}{\partial \nu} \big(L_{\hat{g}} \Phi \big).
 \end{align}
 Using the covariance of the conformal laplacian and the definition of $\Phi$,
 \begin{align} \label{LP1}
 L_{\hat{g}} \Phi = \Gamma_p^{-\frac{n+2}{n-2}} L_g \big( \Gamma_p^{\frac{2}{n-2}} G_p \big).
\end{align}
Let $r(x) = d_g(x,p)$ denote the distance function from $p$ in the metric $g$.  By Lemma 6.4 of \cite{LP}, we can
normalize $\Gamma_p$ so that
\begin{align} \label{GamAs}
\Gamma_p^{\frac{2}{n-2}} = \left\{ \begin{array}{lllll} r^{-2} + O(r) \ \mbox{ if } n = 5, \\
\\
r^{-2} + O(r^2 \log r) \ \mbox{ if } n = 6, \\
\\
r^{-2} + O(r^2) \ \mbox{ if } n = 7.
\end{array}
\right.
\end{align}
Combining this with Proposition \ref{p:green}, for $n = 5,6, 7$ we have
\begin{align} \label{GaGa}
\Gamma_p^{\frac{2}{n-2}} G_p = c_n r^{2-n} + \alpha r^{-2} + O(r^{-1}).
\end{align}
Using Lemma \ref{l:d1d2d4u} and the fact that $R_g = O(r^2)$ in conformal normal coordinates, we get
\begin{align} \label{LGaG} \begin{split}
L_g \big( \Gamma_p^{\frac{2}{n-2}} G_p \big) &= - \Delta_g \big( \Gamma_p^{\frac{2}{n-2}} G_p \big) + \frac{(n-2)}{4(n-1)} R_g \Gamma_p^{\frac{2}{n-2}} G_p \\
&= 2(n-4)\alpha r^{-4} + O(r^{4-n}) \\
&= 2(n-4) \alpha r^{-4} + O(r^{-3}), \ \ \mbox{if }5 \leq n \leq 7.
\end{split}
\end{align}
Note that in dimensions $n \geq 8$ the second term in no longer lower order.
By (\ref{GamAs}),
\begin{align}
\Gamma_p^{-\frac{n+2}{n-2}} = r^{n+2} + O(r^{n+3}),
\end{align}
hence
\begin{align} \label{GGL}
L_{\hat{g}} \Phi = \Gamma_p^{-\frac{n+2}{n-2}} L_g \big( \Gamma_p^{\frac{2}{n-2}} G_p \big) = 2(n-4)\alpha r^{n-2} + O(r^{n-1}).
\end{align}
It is easy to verify that
\begin{align} \label{nur}
\frac{\partial}{\partial \nu} = - \Gamma_p^{-\frac{2}{n-2}} \frac{\partial}{\partial r},
\end{align}
so combining (\ref{GaGa}) and (\ref{GGL}) we find
\begin{align}
\frac{\partial}{\partial \nu} \big(L_{\hat{g}} \Phi \big)\big|_{\partial B_{\delta}} = - 2(n-2)(n-4)\alpha \delta^{n-1} + O(\delta^{n}).
\end{align}
Also, the surface measure transforms by
\begin{align} \label{dS}
\oint_{\partial B_{\delta}} dS_{\hat{g}} = \oint_{\partial B_{\delta}} \Gamma_p^{\frac{2(n-1)}{(n-2)}} dS_g = \omega_{n-1} \delta^{1-n} + O(\delta^{2-n}).
\end{align}
Consequently, the leading boundary term in (\ref{Is1}) is
\begin{align} \label{mainI}
\oint_{\partial B_{\delta}} \frac{\partial}{\partial \nu} \big(\Delta_{\hat{g}} \Phi \big)\ dS_{\hat{g}} = 2(n-2)(n-4)\omega_{n-1} \alpha + o(1).
\end{align}
We can argue as in \cite{HR2} to show that the second boundary integral in (\ref{Is1}) satisfies
\begin{align} \label{obt}
\oint_{\partial B_{\delta}} 4 A_{\hat{g}}(\nabla \Phi, \nu)  dS_{\hat{g}} = o(1),
\end{align}
hence
\begin{align} \label{ppmt}
2(n-2)(n-4)\omega_{n-1} \alpha = (n-4) \int_{M^n \setminus B_{\delta}} |A(\hat{g})|^2 \Phi\ dv_{\hat{g}} + o(1).
\end{align}
It follows that $\alpha \geq 0$.  Moreover, if $\alpha = 0$ then $\hat{g}$ is Ricci-flat, which implies $(X^n,\hat{g})$ is isometric to
flat Euclidean space (see, for example, \cite{Schoen}, Proposition 2, page 492).  This completes the proof.
\end{pf}

%%%%%%%%%%%%%%%%%%%%%%%%%%%%%%%%%%%%%%%%%%%%%%%
\section{The flow} \label{flowSec}
%%%%%%%%%%%%%%%%%%%%%%%%%%%%%%%%%%%%%%%%%%%%%%

\subsection{The initial assumptions} In the following, we assume $(M^n,g_0)$ is a closed Riemannian manifold of dimension $n \geq 5$ with
\begin{align} \label{initpos}  \begin{split}
Q_{g_0} & \mbox{ is semi-positive,} \\
& R_{g_0} \geq 0.
\end{split}
\end{align}
Note that by Lemma \ref{RposLemma}, the assumption on the $Q$-curvature implies $R_{g_0} > 0$.  Also, by Proposition \ref{PPosProp}, $P_{g_0}$ is invertible.  Therefore, we can consider the flow
 \begin{align} \label{uflow}
\left\{ \begin{array}{lll} \displaystyle \frac{\partial u}{\partial t} =  - u + \mu P_{g_0}^{-1} \big( |u|^{\frac{n+4}{n-4}} \big), \\
\\
u(\cdot, 0) = 1,
\end{array}
\right.
\end{align}
where
\begin{align} \nonumber
\mu = \frac{ \int u P_{g_0} u\ dv_0 }{\int |u|^{\frac{2n}{n-4}} dv_0 }.
\end{align}

\begin{lem} \label{STE}
The flow (\ref{uflow}) has a smooth solution for $0 \leq t < T$, where $0 < T \leq \infty$.
\end{lem}

\begin{pf}  Consider the flow
\begin{align} \label{vflow}
\left\{ \begin{array}{lll} \displaystyle \frac{\partial v}{\partial t} =  - v + P_{g_0}^{-1} \big( |v|^{\frac{n+4}{n-4}} \big), \\
\\
v(\cdot, 0) = 1,
\end{array}
\right.
\end{align}
which differs from (\ref{uflow}) by the normalizing term $\mu$.  In fact, these flows just differ by a rescaling in space-time. To see this, suppose $v \in C^{4,\alpha}(M^n \times [0,T))$ is a
solution of (\ref{vflow}), and define
\begin{align} \label{nudef}
\nu = \nu(t) = \displaystyle \frac{ \int v P_{g_0} v\ dv_0 }{\int |v|^{\frac{2n}{n-4}}\ dv_0},
\end{align}
\begin{align} \label{sdef}
s(t) = \int_0^t \nu(\tau)\ d\tau.
\end{align}
Let
\begin{align} \label{urdef}
u(x,t) = e^{s(t) - t}v(x, s(t)).
\end{align}
It is easy to see that $u$ satisfies (\ref{uflow}) on some time interval $[0, \tilde{T})$.

Short-time existence for the flow (\ref{vflow}) follows from the Picard-Lindel\"of theorem on Banach spaces; if we denote $X_{\epsilon} = C^{4,\alpha}(M^n \times [0,\epsilon])$,
then the mapping
\begin{align}
v \mapsto \Psi(v)(x,t) = 1 - \int_0^t v(x,\tau)\ d\tau + \int_0^t P_{g_0}^{-1} (|v|^{\frac{n+4}{n-4}})(x,\tau)\ d\tau
\end{align}
is a contraction  on a small neighborhood of $v_0 \equiv 1$ in $X_{\epsilon}$ for $\epsilon > 0$ small.  A fixed point of $\Psi$ solves (\ref{vflow}).

Note that as (\ref{vflow}) is a non-local ODE in $C^{4,\a}(M)$, there is in general no gain of (spatial) derivatives.
\end{pf}

\begin{proposition} \label{uposProp} For all $0 \leq t < T$,
\begin{align} \label{upos}
u(t,x) > 0.
\end{align}
\end{proposition}

\begin{pf} By (\ref{uflow}),
\begin{align} \label{Pflow} \begin{split}
\frac{\partial}{\partial t}P_{g_0} u &= P_{g_0} (\frac{\partial}{\partial t}u) \\
&= -P_{g_0} u + \mu  |u|^{\frac{n+4}{n-4}},
\end{split}
\end{align}
hence
\begin{align} \label{PODE1}
\frac{\partial}{\partial t}P_{g_0} u \geq -P_{g_0} u.
\end{align}
Integrating this inequality we get
\begin{align} \label{PODE} \begin{split}
P_{g_0} u (t,x) &\geq e^{-t} P_{g_0} u (0,x) \\
&= e^{-t}P_{g_0} (1) \\
&= \frac{n-4}{2} e^{-t} Q_{g_0}(x).
\end{split}
\end{align}
It follows that $P_{g_0} u \geq 0$, and $P_{g_0} u > 0$ somewhere (namely, where the $Q$-curvature is initially positive).
By the strong maximum principle of Theorem \ref{SMP} it follows that $u > 0$ for $t \in [0, T)$.
\end{pf}

\vskip.2in

\begin{remark}  It follows from the proof of Lemma \ref{uposProp} that $Q_g > 0$ for all $t \in (0,T)$: since $u > 0$ for all time, from (\ref{Pflow})
we have
\begin{align*}
\frac{\partial}{\partial t}P_{g_0} u \geq -P_{g_0} u + \mu  u^{\frac{n+4}{n-4}},
\end{align*}
and integrating this we see that $P_{g_0}u > 0$ for $t \in (0,T)$.
\end{remark}

\vskip.2in

%%%%%%%%%%%%%%%%%%%%%%%%%%%%%%%%%%%%%%
\subsection{Variational Properties} \label{SSVar}
%%%%%%%%%%%%%%%%%%%%%%%%%%%%%%%%%%%%%%%

Since $u > 0$ for as long as the flow exists, we can rewrite (\ref{uflow}) as
\begin{align} \label{uposflow}
\frac{\partial}{\partial t}u = - u + \mu P_{g_0} \big( u^{\frac{n+4}{n-4}} \big),
\end{align}
with
\begin{align} \label{mupos}
\mu = \frac{ \int u P_{g_0} u\ dv_0 }{\int u^{\frac{2n}{n-4}}\ dv_0 }.
\end{align}

\begin{lem} \label{Ec}
\begin{align}  \label{dPdt}
\frac{d}{dt} \int u P_{g_0} u\ dv_0 = 0.
\end{align}
\end{lem}

\begin{pf}  From (\ref{uposflow}) and (\ref{mupos}),
\begin{align*}
\frac{d}{dt} \int u P_{g_0} u\ dv_0 &= \int \Big\{ \big(\frac{\partial u}{\partial t}\big) P_{g_0} u + u P_{g_0} \big(\frac{\partial u}{\partial t}\big)\Big\}\ dv_0 \\
&= 2 \int  u P_{g_0} \big(\frac{\partial u}{\partial t}\big)\ dv_0 \\
&= 2 \int u P_{g_0} \Big( - u + \mu P_{g_0}^{-1} \big( u^{\frac{n+4}{n-4}}\big) \Big)\ dv_0 \\
&= 2 \int \Big( - u P_{g_0} u + \mu u P_{g_0} P_{g_0}^{-1} \big( u^{\frac{n+4}{n-4}}\big) \Big)\ dv_0 \\
&= \int \Big\{ - 2 u P_{g_0} u + 2 \mu u^{\frac{2n}{n-4}}\Big\}\ dv_0 \\
&= - 2 \int u P_{g_0} u\ dv_0 + 2 \Big(\frac{ \int u P_{g_0} u\ dv_0 }{\int u^{\frac{2n}{n-4}}\ dv_0 }\Big) \int u^{\frac{2n}{n-4}}\ dv_0 \\
&= 0.
\end{align*}

\end{pf}

To state the next lemma, we denote
\begin{align} \label{fdef}
f = - u + \mu P_{g_0}^{-1} \big( u^{\frac{n+4}{n-4}} \big).
\end{align}
\vskip.2in

\begin{lem} \label{VLem} The conformal volume satisfies
\begin{align} \label{dVdt}
\frac{d}{dt} V = \frac{d}{dt} \int u^{\frac{2n}{n-4}}\ dv_0 = \frac{2n}{n-4} \frac{1}{\mu} \int f P_{g_0} f\ dv_0 \geq 0.
\end{align}
In particular, the volume is increasing along the flow, while $\mu$ and the Paneitz-Sobolev quotient are both decreasing:
\begin{align} \label{mudown} \begin{split}
\frac{d}{dt} \mu &= \frac{d}{dt} \Big(  \frac{ \int u P_{g_0} u\ dv_0}{V} \Big) \leq 0, \\
\frac{d}{dt} \mathcal{F}_{g_0}[u] &= \frac{d}{dt} \Big( \frac{ \int u P_{g_0} u\ dv_0 }{ V^{\frac{n-4}{n}}}\Big) \leq 0.
\end{split}
\end{align}
Finally, the volume is bounded above:
\begin{align} \label{Vabove}
V \leq C_0(g_0).
\end{align}
\end{lem}
\vskip.2in

\begin{pf} To prove the Lemma, we differentiate:
\begin{align} \label{Vp1} \begin{split}
\frac{d}{dt} \int u^{\frac{2n}{n-4}}\ dv_0 &= \frac{2n}{n-4} \int u^{\frac{n+4}{n-4}} \frac{\partial u}{\partial t}\ dv_0 \\
&= \frac{2n}{n-4} \int u^{\frac{n+4}{n-4}} \Big\{ - u + \mu P_{g_0}^{-1}\big( u^{\frac{n+4}{n-4}}\big) \Big\}\ dv_0 \\
&= \frac{2n}{n-4} \int \Big\{ - u^{\frac{2n}{n-4}}+  \mu u ^{\frac{n+4}{n-4}} P_{g_0}^{-1}\big( u^{\frac{n+4}{n-4}}\big) \Big\}\ dv_0.
\end{split}
\end{align}
Note that
\begin{align} \label{fPf} \begin{split}
\int f P_{g_0} f\ dv_0 &= \int \big\{ - u + \mu P_{g_0}^{-1}  \big( u^{\frac{n+4}{n-4}} \big) \big\} \big\{ - P_{g_0} u + \mu u^{\frac{n+4}{n-4}} \big\}\ dv_0 \\
&= \int \Big\{ u P_{g_0} u - \mu u^{\frac{2n}{n-4}} - \mu P_{g_0}^{-1}\big( u^{\frac{n+4}{n-4}} \big) P_{g_0} u  + \mu^2 u ^{\frac{n+4}{n-4}} P_{g_0}^{-1}\big( u^{\frac{n+4}{n-4}}\big)  \Big\}\ dv_0 \\
&= \int \Big\{ - \mu u^{\frac{2n}{n-4}}+  \mu^2 u ^{\frac{n+4}{n-4}} P_{g_0}^{-1}\big( u^{\frac{n+4}{n-4}}\big) \Big\}\ dv_0.
\end{split}
\end{align}
Comparing (\ref{Vp1}) and (\ref{fPf}), we arrive at (\ref{dVdt}).

To see that the volume is bounded above, we use the fact that the Paneitz-Sobolev constant is positive:
\begin{align*}
0 < q_0 \leq \mathcal{F}_{g_0}[u] = V^{-\frac{n-4}{n}} \int u P_{g_0} u\ dv_0 =  V^{-\frac{n-4}{n}} \int u_0 P_{g_0} u_0 \ dv_0 = \frac{n-4}{2} V^{-\frac{n-4}{n}} \int Q_{g_0}\ dv_0,
\end{align*}
hence $V \leq C(g_0)$.

\end{pf}

\begin{corollary} \label{STCor}

We have the space-time estimates
\begin{align} \label{STb} \begin{split}
&\int_0^T  \| f \|_{W^{2,2}}\ dt \leq C_1(g_0), \\
&\int_0^T \Big( \int |f|^{\frac{2n}{n-4}}\ dv_0 \Big)^{\frac{n-4}{n}}\ dt \leq C_2(g_0).
\end{split}
\end{align}
\end{corollary}
\vskip.2in

\begin{pf} From the upper bound on volume we have
\begin{align}
\int_0^T \Big( \int_{M^n} f P_{g_0} f\ dv_0 \Big)\ dt &\leq C_1(g_0).
\end{align}
Since $P_{g_0}$ is positive,
\begin{align*}
\| \phi \|_{W^{2,2}} \approx \int \phi P_{g_0} \phi\ dv_0,
\end{align*}
and the first estimate in (\ref{STb}) follows.  The second estimate follows from the lower bound on the Paneitz-Sobolev quotient.
\end{pf}
\vskip.2in

%%%%%%%%%%
\subsection{Long time existence}
%%%%%%%%%%%%%%%%

\begin{proposition} \label{LTE} The flow (\ref{uflow}) has a smooth solution for all time.  Moreover,
\begin{align}  \label{upperu}
u \leq C' e^{Ct},
\end{align}
where $C, C' > 0$ are constants depending on $g_0$ and the initial datum.
\end{proposition}

\begin{pf} Let $s > 1$.
Since $u > 0$ and $P_{g_0} u > 0$ for as long as the flow exists, by (\ref{Pflow}) we have
\begin{align} \label{intP} \begin{split}
\frac{d}{dt} \int (P_{g_0} u)^s\ dv_0 &= s \int (P_{g_0} u)^{s-1} \frac{\partial}{\partial t} (P_{g_0} u)\ dv_0 \\
&= s \int (P_{g_0} u)^{s-1} \{ - P_{g_0} u + \mu u^{\frac{n+4}{n-4}} \}\ dv_0 \\
&= - s \int (P_{g_0} u)^s + s \mu \int (P_{g_0} u)^{s-1} u^{\frac{n+4}{n-4}}\ dv_0.
\end{split}
\end{align}
For the second integral above we use H\"older's inequality to write
\begin{align} \label{Holder1}
\int (P_{g_0} u)^{s-1} u^{\frac{n+4}{n-4}}\ dv_0 \leq \big( \int (P_{g_0} u)^s\ dv_0 \big)^{\frac{s-1}{s}} \big( \int u^{ \frac{n+4}{n-4}s}\ dv_0 \big)^{\frac{1}{s}}
\end{align}
Assume
\begin{align} \label{sint}
\frac{2n}{n+4} < s < \frac{n}{4}.
\end{align}
Then we can apply H\"older's inequality again to get
\begin{align} \label{holder3}
\big( \int u^{ \frac{n+4}{n-4}s}\ dv_0 \big)^{\frac{1}{s}} \leq \big( \int u^{ \frac{ns}{n-4s}}\ dv_0 \big)^{\frac{n-4s}{ns}} \big( \int u^{ \frac{2n}{n-4} }\ dv_0\big)^{\frac{4}{n}}.
\end{align}
By the Sobolev embedding theorem $W^{s,4} \hookrightarrow L^{ \frac{ns}{n-4s}}$ for $1 < s < n/4$.  Also, since $P_{g_0} > 0$ we have
\begin{align*}
\| u \|_{W^{s,4}} \approx \| P_{g_0} u \|_{L^s}.
\end{align*}
Therefore,
\begin{align} \label{Sob1}
\big( \int u^{ \frac{ns}{n-4s}}\ dv_0 \big)^{\frac{n-4s}{ns}} \leq C_s \big( \int (P_{g_0} u)^s \ dv_0\big)^{\frac{1}{s}}
\end{align}
for $s$ in the range given by (\ref{sint}).  Substituting this into (\ref{holder3}) and using the conformal volume bound of Lemma \ref{VLem} we have

\begin{align} \label{Holder4}
\int (P_{g_0} u)^{s-1} u^{\frac{n+4}{n-4}} \ dv_0\leq C_s  \int (P_{g_0} u)^s\ dv_0.
\end{align}
Substituting this into (\ref{intP}) gives
\begin{align} \label{intP2}
\frac{d}{dt} \int (P_{g_0} u)^s\ dv_0 \leq C_s \int (P_{g_0} u)^s\ dv_0, \hskip.25in \frac{2n}{n+4} < s < \frac{n}{4}.
\end{align}
Integrating this we get
\begin{align} \label{Plp}
\int (P_{g_0} u)^s \ dv_0\leq C_0 e^{C_s t}, \ \ 0 \leq t < T.
\end{align}
By the Sobolev embedding, this implies
\begin{align} \label{uLp1}
\| u \|_{L^{ \frac{ns}{n-4s}}} \leq C_1 e^{C_s' t}.
\end{align}
By choosing $s$ sufficiently close to $n/4$, we conclude that
\begin{align} \label{uLp2}
\| u \|_{L^p} \leq C_3 e^{C_p t},
\end{align}
for any $p > 1$.

Now fix $s > n/4$; say $s = n/4 + 1$.  Returning to (\ref{Holder1}), we have
\begin{align} \label{Holder1a} \begin{split}
\int (P_{g_0} u)^{s-1} u^{\frac{n+4}{n-4}}\ dv_0 &\leq \big( \int (P_{g_0} u\ dv_0)^s \big)^{\frac{s-1}{s}} \big( \int u^{ \frac{n+4}{n-4}s}\ dv_0 \big)^{\frac{1}{s}} \\
&\leq \big( \int (P_{g_0} u)^s \ dv_0\big)^{\frac{s-1}{s}} \big( C_3 e^{C_n t} \big)^{\frac{1}{s}} \\
&\leq C_4 e^{C_5 t} \big( \int (P_{g_0} u)^s \ dv_0\big)^{\frac{s-1}{s}} \\
&\leq \int (P_{g_0} u)^s\ dv_0 + C_6 e^{C_7 t}.
\end{split}
\end{align}
Substituting this into (\ref{intP}) gives
\begin{align} \nonumber
\frac{d}{dt} \int (P_{g_0} u)^s \ dv_0\leq C' e^{C t}, \ \ s = \frac{n}{4} + 1.
\end{align}
Integrating this and using the Sobolev-Kondrakov theorem we conclude
\begin{align} \label{uCa}
\| u \|_{C^{\alpha}} \leq C' e^{C t},
\end{align}
for some $\alpha \in (0,1)$.  This implies \eqref{upperu} and, via (\ref{Pflow}), that the $C^{\alpha}$-norm of $P_{g_0} u$ grows at most exponentially fast. It follows
that $C^{4,\alpha}$-norm of $u$ grows at most exponentially fast, so we cannot have blow-up in finite time.
\end{pf}

%%%%%%%%%%%%%%%%%%%%%%%%%%%%%%%%%
\section{constructing the initial data, part I: $n \geq 8$} \label{Sechighd}
%%%%%%%%%%%%%%%%%%%%%%%%%%%%%%%%%%%

To prove the convergence of the flow we will show that it is possible to construct initial data satisfying the positivity conditions (\ref{initpos}) and with energy below the Euclidean value.  Using a standard argument (see Section \ref{converge}) the latter fact will imply that the flow has a non-zero weak limit which defines a metric of constant $Q$-curvature.

Our first result in
this direction considers the case where the dimension is large (i.e., $n \geq 8$) and the underlying  manifold is not locally conformally flat:

\begin{pro} \label{p:testhighd} Let $(M^n,\bar{g})$  be a closed Riemannian manifold of dimension $n \geq 8$.  Assume  \\

\noindent $(i)$ $Q_{\bar{g}}$ is semi-positive, \\

\noindent $(ii)$ $R_{\bar{g}} \geq 0$, \\

\noindent $(iii)$ $(M^n,\bar{g})$ is not locally conformally flat.  \\

If at $x_0 \in M$ the Weyl tensor $W(x_0)$ is non-zero, then for $\e > 0$ small there exists a function $\psi_{\e} \in C^{\infty}$ and a dimensional constant $c_n$ such that
$$
   \mathcal{F}_{\bar{g}}(\psi_\e) \leq S_n -c_n \e^4 |\log \e \, | \, |W(x_0)|^2 \qquad \hbox{ if } n = 8,
$$
and
$$
  \mathcal{F}_{\bar{g}}(\psi_\e) \leq S_n - c_n \e^4 |W(x_0)|^2 \qquad \hbox{ if } n \geq  9,
$$
where $S_n$ is the Euclidean Paneitz-Sobolev constant:
\begin{align} \label{Sndef}
S_n = \inf_{\varphi \in C_0^{\infty}(\mathbb{R}^n)} \frac{ \int (\Delta_0 \varphi)^2\ dx }{ \Big( \int |\varphi|^{\frac{2n}{n-4}}\ dx \Big)^{\frac{n-4}{n}}}.
\end{align}

Moreover, $\psi_\e$ is positive and induces a conformal metric $h = \psi_\e^{\frac{4}{n-4}}\bar{g}$ with the following properties: \\

\noindent $(i')$ $Q_{h}$ is semi-positive,  \\

\noindent $(ii')$ $R_{h} > 0$, \\

\noindent $(iii')$
\begin{align} \label{Fh1} \begin{split}
   \mathcal{F}_{h}(1) &\leq S_n -c_n \e^4 |\log \e \, | \, |W(x_0)|^2 \ \ \hbox{ if } n = 8, \\
  \mathcal{F}_{h}(1) &\leq S_n - c_n \e^4 |W(x_0)|^2 \qquad \ \ \ \ \hbox{     if } n \geq  9.
  \end{split}
  \end{align}
\end{pro}

\

\

\begin{pf}  Let $\tilde{g} = \varphi^{\frac{4}{n-4}} \bar{g}$ denote the metric
satisfying the conformal normal coordinate conditions of  \cite{LP} at $x_0$ (we assume $\varphi$ is globally defined). Consider the test function
in Section 6 of \cite{ER} defined by
$$
  \tilde{u}_\e(x) = \frac{\eta(x) \varphi(x)}{\left( \e ^2 + d_{\tilde{g}}(x,x_0)^2 \right)^{\frac{n-4}{2}}},
$$
where  $\eta(x)$ is a cut-off function with support in a ball
$B_{2\d}(x_0)$, identically equal to 1 in $B_\d(x_0)$.

In Section 7 of \cite{ER} it was shown that, for $\e > 0$ small one has the estimates
$$
   \mathcal{F}_{\bar{g}}(\tilde{u}_\e) \leq S_n - C(n) \e^4 |\log \e \, | \, |W(x_0)|^2 \qquad \hbox{ if } n = 8
$$
and
$$
  \mathcal{F}_{\bar{g}}(\tilde{u}_\e) \leq S_n - C(n) \e^4 |W(x_0)|^2 \qquad \hbox{ if } n \geq  9.
$$
We will show that it is possible to modify these test functions in order to produce a strictly positive conformal factor
which defines a metric with semi-positive $Q$ and positive scalar curvatures, while preserving the property of the Paneitz-Sobolev quotient being below the Euclidean value.  We begin with the following lemma:

\begin{lem}\label{l:panue} If $\tilde{g}$ is as above, if we set
$$
  {u}_\e(x) = \frac{\eta(x)}{\left( \e ^2 + d_{\tilde{g}}(x,x_0)^2 \right)^{\frac{n-4}{2}}},
$$
then
\begin{equation}\label{eq:ptildeg}
   P_{\tilde{g}}(u_\e) = \frac{ n(n-4)(n^2-4)  \e^4}{\left( \e ^2 + |x|^2
   \right)^{\frac{n+4}{2}}} + \frac{O(1)}{\left( \e^2 + r^2
   \right)^{\frac{n-4}{2}}} \qquad \quad \hbox{ in } B_{2\d}(x_0).
\end{equation}
\end{lem}

\noindent Notice that, by the conformal covariance of the Paneitz operator we have that $\mathcal{F}_{\bar{g}}(\tilde{u}_\e) =
\mathcal{F}_{\tilde{g}}(u_\e)$. From now on we will work in the metric $\tilde{g}$.

\

\begin{pf} The estimate is trivial in $B_{2\d}(x_0) \setminus B_\d(x_0)$ (where the second term
in the r.h.s. of \eqref{eq:ptildeg} dominates the first one). It is therefore
sufficient to prove it in $B_\d(x_0)$, where $\eta$ is identically equal to $1$,
and hence here it is enough to estimate
$$
  P_{\tilde{g}} \left( \left( \e ^2 + r^2 \right)^{\frac{4-n}{2}} \right).
$$
Let us first consider the bi-Laplacian term: for a radial function $f(r)$ in conformal normal
coordinates we have that
$$
  \D f(r) = \frac{1}{\sqrt{\det g}} \pa_r \left( \sqrt{\det g} \, \pa _r f \right) =
  f'' + \frac{n-1}{r} f' + O(r^{N-1}) f',
$$
where $N \geq 5$. Therefore, if $\Delta_0$ denotes the Euclidean Laplacian
$$
  \D^2 f(r) = \D ^2_{0} f + O(r^{N-1}) f''' + O(r^{N-2}) f'' + O(r^{N-3}) f'.
$$
By an explicit computation we find that, if $f(r) = \left( \e ^2 + r^2 \right)^{\frac{4-n}{2}}$, then
$$
  \D ^2_{0} f  = n(n-4)(n^2-4) \frac{\e ^4}{\left( \e ^2 + r^2 \right)^{\frac{n+4}{2}}} :=
  b_n\frac{\e ^4}{\left( \e ^2 + r^2 \right)^{\frac{n+4}{2}}},
$$
and (for a dimensional constant $a_n$)
$$
  |f'| \leq   \frac{a_n \, r}{\left( \e ^2 + r^2 \right)^{\frac{n-2}{2}}}; \qquad
  |f''| \leq   \frac{a_n}{\left( \e ^2 + r^2 \right)^{\frac{n-2}{2}}}; \qquad
    |f'''| \leq \frac{a_n \, r}{\left( \e ^2 + r^2 \right)^{\frac{n}{2}}}.
$$
Therefore we obtain that
$$
  \D^2 f(r) = b_n \frac{\e ^4}{\left( \e ^2 + r^2 \right)^{\frac{n+4}{2}}}
  + \frac{O(r^{N-2})}{\left( \e ^2 + r^2 \right)^{\frac{n-2}{2}}}
  = \frac{\left( b_n \e ^4
    + O \left( r^{N-2} (\e^2 + r^2)^3  \right) \right)}{\left( \e ^2 + r^2 \right)^{\frac{n+4}{2}}}.
$$
Next, we check the lower order terms of the Paneitz operator.  Recall
\begin{align} \label{Pforms}
P f = \Delta^2 f + c_1 R_{ij} \nabla_i \nabla_j f + c_2 R \Delta f + c_3 \langle \nabla R , \nabla f \rangle + c_4 Q f,
\end{align}
where $R_{ij}$ are the components of the Ricci tensor, and the $c_i$'s are dimensional constants. In conformal normal coordinates,
\begin{align} \label{curvCN}
Ric(\frac{\partial}{\partial r},\frac{\partial}{\partial r}) = O(r^2), \ \ R = O(r^2), \ \ |\nabla R| = O(r), \ \ |Q| = O(1).
\end{align}
Therefore, the terms in $Pf$ involving first and second derivatives of $f$ are of the order
$$
  r f' + r^2 f'',
$$
which are bounded by
$$
  \frac{O(r^2)}{\left( \e ^2 + r^2 \right)^{\frac{n-2}{2}}} =
  \frac{O(1)}{\left( \e ^2 + r^2 \right)^{\frac{n-4}{2}}}.
$$
The term $Q_{\tilde{g}} u_\e$ is bounded  by a  constant times $f$, namely
$$
  \frac{O(1)}{\left( \e ^2 + r^2 \right)^{\frac{n-4}{2}}}.
$$
In conclusion we find that
$$
  P_{\tilde{g}} \left( u_\e \right) = \frac{ b_n \,  \e^4}{\left( \e ^2 + |x|^2 \right)^{\frac{n+4}{2}}} +
   \frac{O \left( r^{N-2} (\e^2 + r^2)^3  \right)}{\left( \e ^2 + r^2 \right)^{\frac{n+4}{2}}}
   +   \frac{O(1)}{\left( \e^2 + r^2 \right)^{\frac{n-4}{2}}}.
$$
For $N$ sufficiently large the second term in the r.h.s. can be absorbed into the third, so we obtain
the  desired estimate.
\end{pf}

\

\noindent Recalling the invertibility of $P$ from Proposition \ref{PPosProp}, we consider next the function $\hat{u}_\e$ defined by the equation
\begin{equation}\label{eq:hatue}
  P_{\tilde{g}} \hat{u}_\e = \eta(x) \frac{b_n \, \e^4}{\left( \e ^2 + |x|^2 \right)^{\frac{n+4}{2}}}.
\end{equation}
We aim to estimate the difference between this new function and $u_\e$.

\begin{lem}\label{l:diff}
If $\hat{u}_\e$ is as above, let us set
$$
 v_\e = \hat{u}_\e - u_\e .
$$
Then there exists $C > 0$ such that in $B_{2\d}(x_0)$ we have the estimates
$$
 | v_\e | \leq C \left( \e^2 + |x|^2 \right)^{\frac{8-n}{2}} \qquad
 \quad \hbox{ if } n > 8;
$$
$$
    | v_\e | \leq C \log \left( \frac{1}{\e^2 + |x|^2} \right) \qquad
     \quad \hbox{ if } n = 8.
$$
On $M \setminus B_{2 \d}(x_0)$ we have simply
$$
  |v_\e | \leq C.
$$
\end{lem}

\begin{pf}
We notice that, by Lemma \ref{l:panue}
\begin{equation}\label{eq:pve}
 P_{\tilde{g}}(v_\e) = P_{\tilde{g}} \left( \hat{u}_\e - u_\e \right) = \eta(x) \frac{b_n \,
  \e^4}{\left( \e ^2 + |x|^2 \right)^{\frac{n+4}{2}}}   - P_{\tilde{g}} u_\e
  = \frac{O(1)}{\left( \e^2 + r^2 \right)^{\frac{n-4}{2}}}.
\end{equation}
Recall also that the r.h.s. is supported in $B_{2\d}(x_0)$ as $\eta$ and $u_\e$ are. We estimate now
the convolution of the r.h.s. with the Green's function of the Paneitz operator, which
is bounded above by $O(1)/d_{\tilde{g}}(x,y)^{n-4}$.

For $n = 8$ we can divide between the regime $|x| = O(\e)$ and $|x| \geq C_0 \e$ for a large constant $C_0$. When $n = 8$ and
$|x| = O(\e)$ the convolution is bounded by
$$
  C \int_{|y| \leq 1} \frac{1}{|x-y|^4} \frac{dy}{\left( \e^2 + |y|^2 \right)^{\frac{n-4}{2}}}.
$$
By a change of variables ($y = \e \, w)$ one finds that this integral can be controlled by
$$
  C \int_{|w| \leq 1/\e} \frac{1}{|\overline{x}-w|^4} \frac{dw}{\left( 1 + |w|^2 \right)^{\frac{n-4}{2}}},
$$
where $|\overline{x}| = O(1)$. One can easily see that the latter integral is of order
$\log \frac{1}{\e}$. On the other hand, for $|x| \geq C_0 \e$ we can write that
$$
  \int_{|y| \leq 1} \frac{1}{|x-y|^4} \frac{dy}{\left( \e^2 + |y|^2 \right)^{\frac{n-4}{2}}}
  = \int_{|y| \leq \frac{1}{|x|}} \frac{dw}{\left| \frac{x}{|x|} - w \right|^4
  \left( \frac{\e^2}{|x|^2} + |w|^2 \right)^2} \leq \log \frac{1}{|x|}.
$$
In conclusion for $n = 8$ we get
$$
  |v_\e|(x) \leq \log \left( \frac{1}{\e^2 + |x|^2} \right); \qquad \quad x \in
  B_{2 \d}(x_0).
$$
The estimate on $v_\e$ outside $B_{2 \d}(x_0)$ is immediate.

\

Let us consider now the case $n \geq 9$. We distinguish again between $|x| = O(\e)$
and $|x| \geq C_0 \e$. In the former case we get, similarly to before
$$
  C \int_{|y| \leq 1} \frac{1}{|x-y|^{n-4}} \frac{dy}{\left( \e^2 + |y|^2 \right)^{\frac{n-4}{2}}} =
  C \e^{8-n} \int_{|w| \leq 1/\e} \frac{1}{|\overline{x}-w|^{n-4}} \frac{dw}{\left( 1 + |w|^2 \right)^{\frac{n-4}{2}}},
$$
with $|\overline{x}| = O(1)$. The last integral is uniformly bounded for $n > 9$.

If the case $|x| \geq C_0 \e$ we write
$$
  \int_{|y| \leq 1} \frac{1}{|x-y|^{n-4}} \frac{dy}{\left( \e^2 + |y|^2 \right)^{\frac{n-4}{2}}}
  = |x|^{8-n} \int_{|y| \leq \frac{1}{|x|}} \frac{dw}{\left| \frac{x}{|x|} - w \right|^{n-4}
  \left( \frac{\e^2}{|x|^2} + |w|^2 \right)^{\frac{n-4}{2}}} \leq C |x|^{8-n}.
$$
In conclusion for $n > 8$ we get
$$
  |v_\e|(x) \leq C \left( \e^2 + |x|^2 \right)^{\frac{8-n}{2}};  \qquad \quad x \in
  B_{2 \d}(x_0).
$$
The estimate on $v_\e$ outside $B_{2 \d}(x_0)$ is again quite easy.

This concludes the proof.
\end{pf}

\

\noindent We check next the effect of the correction $v_\e$ on the Paneitz-Sobolev quotient, and in particular how
much it deviates from the Euclidean one.

\begin{lem}  \label{deviation}
One has that
$$
  \mathcal{F}_{\tilde{g}}(\hat{u}_\e) = \mathcal{F}_{\tilde{g}}(u_\e) + o(\e^4 |\log \e \, |) \qquad \quad \hbox{ for } n = 8,
$$
and
$$
  \mathcal{F}_{\tilde{g}}(\hat{u}_\e) = \mathcal{F}_{\tilde{g}}(u_\e) + o(\e^4) \qquad \quad \hbox{ for } n \geq 9.
$$
\end{lem}

\begin{pf}
Calling $\mathcal{N}$ and $\mathcal{D}$ the numerator and the
denominator in the quotient, we have that
$$
  \mathcal{N}(\hat{u}_\e) = \int_M \hat{u}_\e P_{\tilde{g}} \hat{u}_\e dv_{\tilde{g}} =
  \int_M {u}_\e P_{\tilde{g}} {u}_\e dv_{\tilde{g}} + 2 \int_M v_\e P_{\tilde{g}} u_\e
  dv_{\tilde{g}} + \int_M v_\e P_{\tilde{g}} v_\e dv_{\tilde{g}}.
$$
The second term by Lemma \ref{l:panue} can be estimated as
$$
   2 \int_M v_\e \left( \frac{b_n \e^4}{\left( \e^2 + |x|^2 \right)^{\frac{n+4}{2}}} +
   \frac{O(1)}{\left( \e^2 + |x|^2 \right)^{\frac{n-4}{2}}}\right)
     dv_{\tilde{g}}.
$$
By Lemma \ref{l:diff} we can write that
$$
  \int_M v_\e \frac{O(1)}{\left( \e^2 + |x|^2 \right)^{\frac{n-4}{2}}} dv_{\tilde{g}}
  \leq C \int_{B_1(0)} \log \left( \frac{1}{\e^2 + |x|^2} \right)
  \frac{dx}{\left( \e^2 + |x|^2 \right)^{\frac{n-4}{2}}}
$$
for $n = 8$, and
\begin{equation}\label{eq:lhs2}
\int_M v_\e \frac{O(1)}{\left( \e^2 + |x|^2 \right)^{\frac{n-4}{2}}} dv_{\tilde{g}}
    \leq C \int_{B_1(0)}
    \frac{dx}{\left( \e^2 + |x|^2 \right)^{\frac{n-4}{2}+\frac{n-8}{2}}}
\end{equation}
for $n \geq 9$. In the former case, using the change of variables $s = \e^2 + |x|^2$
we can write that
\begin{eqnarray}\label{eq:sim} \nonumber
 \int_{B_1(0)} \log \left( \frac{1}{\e^2 + |x|^2} \right)
    \frac{dx}{\left( \e^2 + |x|^2 \right)^{\frac{8-4}{2}}} & \leq  &
    C \int_{0}^1 \log \left( \frac{1}{\e^2 + |x|^2} \right)
       (\e^2 + |x|^2)^3 |x| \frac{d|x|}{\left( \e^2 + |x|^2 \right)^{2}} \\
       & \leq & C \int_{0}^1 \log \left( \frac{1}{s} \right)
              s \, ds \leq C.
\end{eqnarray}
In the latter case, one can also easily check boundedness of the l.h.s. of \eqref{eq:lhs2} using a change of variables.
In either case we can write that
$$
  2 \int_M v_\e P_{\tilde{g}} u_\e
    dv_{\tilde{g}}  = 2 \int_M v_\e \frac{b_n \e^4}{\left( \e^2 +
    |x|^2 \right)^{\frac{n+4}{2}}} dv_{\tilde{g}}  +  O(1).
$$
In conclusion we get that
$$
\mathcal{N}(\hat{u}_\e)  = \int_M  {u}_\e P_{\tilde{g}} {u}_\e  dv_{\tilde{g}} + 2 b_n \, \e^4
\int_M u_\e^{\frac{n+4}{n-4}} v_\e dv_{\tilde{g}} + O(1).
$$

We turn next to  the denominator $\mathcal{D}$, for which we have
$$
  \mathcal{D}(\hat{u}_\e) = \left( \int_M |u_\e + v_\e|^{\frac{2n}{n-4}}
  dv_{\tilde{g}} \right)^{\frac{n-4}{n}}.
$$
In $B_\d(x_0)$, by Lemma \ref{l:diff} and the explicit expression of $u_\e$, we have that
$|v_\e| \leq C |u_\e|$, so a Taylor expansion gives that
$$
  \left| |u_\e + v_\e|^{\frac{2n}{n-4}}  - u_\e^{\frac{2n}{n-4}} -
   \frac{2n}{n-4}  u_\e^{\frac{n+4}{n-4}}  v_\e \right|
  \leq C u_\e^{\frac{8}{n-4}}  v_\e^2  \qquad \qquad \hbox{ in } B_\d(x_0).
$$
Hence, using again Lemma \ref{l:diff} and the explicit expression of $u_\e$ we can write that
\begin{eqnarray*}
 \int_M |u_\e + v_\e|^{\frac{2n}{n-4}}
   dv_{\tilde{g}} & = & \int_{B_\d(x_0)} |u_\e + v_\e|^{\frac{2n}{n-4}}
     dv_{\tilde{g}} + \int_{M \setminus B_\d(x_0)} |u_\e + v_\e|^{\frac{2n}{n-4}}
       dv_{\tilde{g}} \\ & = & \int_{B_\d(x_0)} \left( u_\e^{\frac{2n}{n-4}} +
          \frac{2n}{n-4}  u_\e^{\frac{n+4}{n-4}}  v_\e + O(u_\e^{\frac{8}{n-4}}  v_\e^2 ) \right)
            dv_{\tilde{g}} + O(1) \\ & = & \int_{M} u_\e^{\frac{2n}{n-4}} dv_{\tilde{g}}
        + \frac{2n}{n-4} \int_M  u_\e^{\frac{n+4}{n-4}}  v_\e     dv_{\tilde{g}}
     + \int_{B_\d(x_0)}  O(u_\e^{\frac{8}{n-4}}  v_\e^2 ) dv_{\tilde{g}} + O(1).
\end{eqnarray*}
Similarly to \eqref{eq:sim} for $n = 8$ and with a change of variables for $n \geq 9$ we
obtain
$$
  \int_{B_\d(x_0)}  O(u_\e^{\frac{8}{n-4}}  v_\e^2 ) dv_{\tilde{g}} = O(1),
$$
and hence we find
$$
  \mathcal{D}(\hat{u}_\e) = \left( \int_{M} u_\e^{\frac{2n}{n-4}} dv_{\tilde{g}}
          + \frac{2n}{n-4} \int_M  u_\e^{\frac{n+4}{n-4}}  v_\e     dv_{\tilde{g}}
       +  O(1) \right)^{\frac{n-4}{n}}.
$$
In conclusion we deduce
$$
  \mathcal{F}_{\tilde{g}}(\hat{u}_\e) = \frac{\int_M  {u}_\e P_{\tilde{g}} {u}_\e  dv_{\tilde{g}} + 2 b_n \, \e^4
  \int_M u_\e^{\frac{n+4}{n-4}} v_\e dv_{\tilde{g}} + O(1)}{\left( \int_{M} u_\e^{\frac{2n}{n-4}} dv_{\tilde{g}}
            + \frac{2n}{n-4} \int_M  u_\e^{\frac{n+4}{n-4}}  v_\e     dv_{\tilde{g}}
         +  O(1) \right)^{\frac{n-4}{n}}},
$$
which means
$$
  \mathcal{F}_{\tilde{g}}(\hat{u}_\e) = \frac{\mathcal{N}(u_\e)}{\mathcal{D}(u_\e)} \frac{1 +
  \frac{2 b_n \, \e^4
    \int_M u_\e^{\frac{n+4}{n-4}} v_\e dv_{\tilde{g}}}{\int_M  {u}_\e P_{\tilde{g}} {u}_\e  dv_{\tilde{g}}}
    + \frac{O(1)}{\int_M  {u}_\e P_{\tilde{g}} {u}_\e  dv_{\tilde{g}}}}{\left( 1
                + \frac{2n}{n-4} \frac{\int_M  u_\e^{\frac{n+4}{n-4}}  v_\e     dv_{\tilde{g}}}{ \int_{M} u_\e^{\frac{2n}{n-4}} dv_{\tilde{g}} }
             +  \frac{O(1)}{ \int_{M} u_\e^{\frac{2n}{n-4}} dv_{\tilde{g}} } \right)^{\frac{n-4}{n}}}.
$$
Notice by Lemma \ref{l:panue},
$$
  \int_M u_\e P_{\tilde{g}} u_\e dv_{\tilde{g}} = b_n \e^4 (1 + o_\e(1)) \int_M u_\e^{\frac{2n}{n-4}} dv_{\tilde{g}},
$$
which implies
\begin{equation}\label{eq:dddd}
\mathcal{F}_{\tilde{g}}(\hat{u}_\e) = \frac{\mathcal{N}(u_\e)}{\mathcal{D}(u_\e)} \frac{1 +
  \frac{2 \int_M u_\e^{\frac{n+4}{n-4}} v_\e dv_{\tilde{g}}}{ (1 + o_\e(1)) \int_M u_\e^{\frac{2n}{n-4}} dv_{\tilde{g}}}
    + \frac{O(1)}{\int_M u_\e^{\frac{2n}{n-4}} dv_{\tilde{g}}}}{\left( 1
                + \frac{2n}{n-4} \frac{\int_M  u_\e^{\frac{n+4}{n-4}}  v_\e     dv_{\tilde{g}}}{ \int_{M} u_\e^{\frac{2n}{n-4}} dv_{\tilde{g}} }
             +  \frac{O(1)}{ \int_{M} u_\e^{\frac{2n}{n-4}} dv_{\tilde{g}} } \right)^{\frac{n-4}{n}}}.
\end{equation}
Now notice that we have the asymptotics
$$
  \int_M u_\e^{\frac{2n}{n-4}} dv_{\tilde{g}} \simeq \e^{-n} , \qquad \hbox{ and } \qquad
  \int_M u_\e^{\frac{n+4}{n-4}} v_\e dv_{\tilde{g}} =
    \begin{cases}
    O(\e^{-4} | \log \e \, |)  &    \hbox{ for } n = 8 \\
     O(\e^{4-n}) &
      \hbox{ for } n \geq 9.
    \end{cases}
$$
Hence from a Taylor expansion of the denominator in \eqref{eq:dddd} we find that
$$
  \mathcal{F}_{\tilde{g}}(\hat{u}_\e)  = \begin{cases}
      (1 + o(\e^{4} | \log \e \, |)) \, \mathcal{F}_{\tilde{g}}(u_{\e}) &    \hbox{ for } n = 8, \\
       (1 + o(\e^{4})) \, \mathcal{F}_{\tilde{g}}(u_{\e}) &
        \hbox{ for } n \geq 9.
      \end{cases}
$$
This concludes the proof.
\end{pf}

\
\begin{lem}\label{u:pos}  $\hat{u}_\e$ is positive.
\end{lem}

\begin{pf}
By the defining equation for $\hat{u}_\e$ and the conformal covariance of the Paneitz operator,
\begin{align*}
P_{\bar{g}}(\varphi \, \hat{u}_\e) &= \varphi^{\frac{n+4}{n-4}} P_{\tilde{g}} (\hat{u}_\e ) \\
&= \varphi(x)^{\frac{n+4}{n-4}} \eta(x) \frac{n(n-4)(n^2-4) \e^4}{\left( \e ^2 + |x|^2 \right)^{\frac{n+4}{2}}} \\
&\geq 0.
\end{align*}
Since $Q_{\bar{g}}\geq 0$ and $R_{\bar{g}} > 0$, by the strong maximum principle of Theorem \ref{SMP} it follows that $\hat{u}_\e > 0$.
\end{pf}

\

Let
\begin{align} \label{psidef}
\psi_{\e} = \varphi \, \hat{u}_\e,
\end{align}
and
\begin{align} \label{hdef1}
h = \psi_{\e}^{4/(n-4)} \bar{g} = \hat{u}_\e^{4/(n-4)} \tilde{g}.
\end{align}

\

\begin{lem}\label{p:pos}  The scalar curvature of the metric $h$ is positive.
\end{lem}

\begin{pf}
For $0 \leq s \leq 1$ let
\begin{align} \label{wt}
w_s = (1-s)\varphi^{-1} + s  \hat{u}_\e,
\end{align}
and
\begin{align} \label{gs}
h_s = w_s^{4/(n-4)}\tilde{g}.
\end{align}
Then
\begin{align*}
h_0 &= \varphi^{-\frac{4}{n-4}}\tilde{g} \\
&= \varphi^{-\frac{4}{n-4}}\big\{ \varphi^{\frac{4}{n-4}}g \big\} \\
&= \bar{g},
\end{align*}
and $h_1 = h$.  Observe that the $Q$-curvature of $h_s$ is semi-positive; this follows from the fact that
\begin{align} \label{Qpos} \begin{split}
P_{\tilde{g}}(w_s) &= (1 - s) P_{\tilde{g}}(\varphi^{-1}) + s P_{\tilde{g}} \hat{u}_\e \\
&= (1-s) P_{\varphi^{\frac{4}{n-4}}\bar{g}} (\varphi^{-1}) + s P_{\tilde{g}} \hat{u}_\e \\
&= (1-s) \varphi^{-\frac{n+4}{n-4}} P_{\bar{g}} (1) + s P_{\tilde{g}} \hat{u}_\e  \\
&= (1-s)\frac{n-4}{2} \varphi^{-\frac{n+4}{n-4}} Q_{\bar{g}} + s P_{\tilde{g}} \hat{u}_\e  \\
&\geq 0,
\end{split}
\end{align}
and clearly $P_{\tilde{g}}(w_s) > 0$ somewhere. Also, note that $R_{h_0} = R_{\bar{g}} > 0$.  Therefore, if there is a $s_1 \in (0,1]$ such that
\begin{align} \label{minR}
\min R_{h_{s_0}} = 0,
\end{align}
then this would contradict Lemma \ref{RposLemma}.  It follows that $R_h > 0$.
\end{pf}

\

To conclude the proof of Proposition \ref{p:testhighd}, we point out that the defining equation for $\hat{u}_\e$ clearly
implies that $P_{\tilde{g}} \hat{u}_\e \geq 0$, with $P_{\tilde{g}} \hat{u}_\e > 0$ near $x_0$.  In particular, the
$Q$-curvature of $h$ is non-negative everywhere, and positive near $x_0$.  We conclude that $(i')$ and $(ii')$ both hold.
Finally, that (\ref{Fh1}) holds follows from Lemma \ref{deviation} and the conformal invariance of $\mathcal{F}$.
\end{pf}

%%%%%%%%%%%%%%%%%%%%%%%%%%%%%%%%%%%%%%%%%%%%%%%%%%%%%%%%%%%%%
\section{constructing the initial data, part II: $n = 5,6,7$ or $\overline{g}$
locally conformally flat}  \label{Seclowd}
%%%%%%%%%%%%%%%%%%%%%%%%%%%%%%%%%%%%%%%%%%%%%%%%%%%%%%%%%%%%%%%%

 In low dimensions (i.e., $n = 5, 6, 7$ or in the locally conformally flat case) the Green's function plays a role in the Paneitz-Sobolev quotient expansion, just
 as for Yamabe's problem in Schoen's work \cite{Schoen}. Using Theorem \ref{PMT}, we will prove

\begin{pro} \label{p:testlowd} Let $(M^n,\bar{g})$  be a closed Riemannian manifold of dimension $n$, with $n= 5, 6$, or $7$; or let $(M^n, \overline{g})$ be locally conformally flat of dimension $n \geq 5$.  Assume  \\

\noindent $(i)$ $Q_{\bar{g}}$ is semi-positive, \\

\noindent $(ii)$ $R_{\bar{g}} \geq 0$.  \\

If $(M^n, \ov{g})$ is not conformally equivalent to the standard sphere, then for $\epsilon > 0$ small and every $x_0 \in M$, there exists a function $\psi_{\e} \in C^{\infty}$ and a  constant $c_{x_0} > 0$ such that
\begin{align} \label{Fbarp2}
   \mathcal{F}_{\bar{g}}(\psi_\e) \leq S_n - c_{x_0} \e^{n-4}.
\end{align}

Moreover, $\psi_\e$ is positive and induces a conformal metric $h = \psi_\e^{\frac{4}{n-4}}\bar{g}$ with the following properties: \\

\noindent $(i')$ $Q_{h}$ is semi-positive,  \\

\noindent $(ii')$ $R_{h} > 0$, \\

\noindent $(iii')$
$$
   \mathcal{F}_{h}(1) \leq S_n - c_{x_0} \e^{n-4}.
$$

\end{pro}

\

\

\begin{pf}
If $n = 5, 6$ or $7$, we let $\varphi$ be as in the proof of Proposition \ref{p:testhighd}. If $\overline{g}$
is locally conformally flat, we choose $\varphi$ so that $\tilde{g} = \varphi^{\frac{4}{n-4}} \overline{g}$
is flat near $x_0$.  We still consider the functions $\hat{u}_\e$ as in \eqref{eq:hatue}, with base point $x_0$, and we try to deduce estimates by
evaluating the Paneitz operator on an approximation.

We consider a cut-off function $\tilde{\chi}_{\tilde{\d}}(x) = \tilde{\chi}(x/\tilde{\d})$, where $\tilde{\chi}$ is a cut-off
function equal to $1$ in $B_1$ and equal to zero outside $B_2$. We then
construct an approximate solution $\check{u}_\e$ defining it by
$$
 \check{u}_\e := \tilde{\chi}_{\tilde{\d}}(u_\e + \beta) + (1 - \tilde{\chi}_{\tilde{\d}}) \ov{G}_{x_0},
$$
where $\beta = \beta_{x_0} = \frac{1}{c_n} \a_{x_0} > 0$, $\a_{x_0}$ appears in the expansion of $G_{x_0}$ in \eqref{Gexp}, and
$\ov{G}_{x_0} = \frac{1}{c_n} G_{x_0}$ with $\tilde{\d} \ll \d$.
Notice that by the positivity of the Green's function (see Section 2), the function $\check{u}_\e$
is positive on $M$. The Paneitz operator on $u_\e$ was
already estimated in the previous Section. We have the following result
concerning an estimate of $P_{\tilde{g}} \check{u}_\e$ in $B_{2\tilde{\d}} \setminus B_{\tilde{\d}}$:

\begin{lem}\label{l:Pgd}
There exists a constant $C > 0$ such that
$$
  |P_{\tilde{g}} \check{u}_\e| \leq C \tilde{\d}^{-3} \qquad \hbox{ in } B_{2\tilde{\d}} \setminus B_{\tilde{\d}}.
$$
\end{lem}

\begin{pf}
We can write
$$
  \check{u}_\e = \ov{G}_{x_0} + \tilde{\chi}_{\tilde{\d}} \left( u_\e + \beta - \ov{G}_{x_0} \right),
$$
and hence it follows that, in $B_{2 \tilde{\d}} \setminus B_{\tilde{\d}}$
\begin{eqnarray*}
|P_{\tilde{g}} \check{u}_\e| & \leq & |\n^4 \tilde{\chi}_{\tilde{\d}}| \left| u_\e + \beta - \ov{G}_{x_0}  \right|
  + |\n^3 \tilde{\chi}_{\tilde{\d}}| \left| \nabla \left( u_\e + \beta - \ov{G}_{x_0} \right) \right|
  \\ & + & |\n^2 \tilde{\chi}_{\tilde{\d}}| \left| \nabla^2 \left( u_\e + \beta - \ov{G}_{x_0} \right) \right| \\ & + &
  |\n \tilde{\chi}_{\tilde{\d}}| \left| \nabla^3 \left( u_\e + \beta - \ov{G}_{x_0} \right) \right| +
  \left| P_{\tilde{g}} \left( u_\e + \beta - \ov{G}_{x_0} \right) \right|.
\end{eqnarray*}
As $\chi_{\tilde{\d}}$ satisfies the estimates
$$
  |\n \tilde{\chi}_{\tilde{\d}}| \leq \frac{C}{\tilde{\d}}; \qquad |\n^2 \tilde{\chi}_{\tilde{\d}}| \leq \frac{C}{\tilde{\d}^2};
  \qquad |\n^3 \tilde{\chi}_{\tilde{\d}}| \leq \frac{C}{\tilde{\d}^3}; \qquad |\n^4 \tilde{\chi}_{\tilde{\d}}| \leq \frac{C}{\tilde{\d}^4},
$$
it will be sufficient to show that in $B_{2\tilde{\d}} \setminus B_{\tilde{\d}}$
$$
    \left| u_\e + \beta - \ov{G}_{x_0}  \right| \leq C
    \tilde{\d}; \qquad
    \left| \nabla \left( u_\e + \beta - \ov{G}_{x_0} \right) \right| \leq C;
    \qquad \left| \nabla^2 \left( u_\e + \beta - \ov{G}_{x_0} \right) \right| \leq \frac{C}{\tilde{\d}};
$$
$$
    \left| \nabla^3 \left( u_\e + \beta - \ov{G}_{x_0} \right) \right| \leq \frac{C}{\tilde{\d}^2};
    \qquad  \left| P_{\tilde{g}} \left( u_\e + \beta - \ov{G}_{x_0} \right) \right| \leq \frac{C}{\tilde{\d}^3}.
$$
We begin with the last inequality: we have that
$$
  P_{\tilde{g}} \left( u_\e + \beta - \ov{G}_{x_0} \right) = P_{\tilde{g}} \left( u_\e + \beta \right)
  = P_{\tilde{g}} u_\e + \beta\, Q_{\tilde{g}} = O(\tilde{\d}^{-3}),
$$
where in the above formula we used Lemma \ref{l:panue} (in the locally conformally
flat case it is obvious).
To prove the remaining estimates we use the fact that in $B_{2\tilde{\d}}$
$$
  u_\e + \beta - \ov{G}_{x_0} = \left( \e^2 + |x|^2 \right)^{\frac{4-n}{2}} - |x|^{4-n} + O_p(1),
$$
by Proposition \ref{p:green}. We remark that in the locally conformally flat case the
above estimate simply follows from the fact that, in the metric $\tilde{g}$,
$\ov{G}_{x_0}(\cdot) - \beta - d_{\tilde{g}}(x_0, \cdot)$ is a smooth bi-harmonic function.

{From} a Taylor expansion of  $\left( \e^2 + |x|^2 \right)^{\frac{4-n}{2}}$ one easily finds that
$$
  u_\e + \beta - \ov{G}_{x_0} = O(\e^2 |x|^{2-n}) + O_p(1).
$$
This implies the conclusion.
\end{pf}

\

Combining the estimates of Lemmas \ref{l:panue} and \ref{l:Pgd} we find
$$
 \left| P_{\tilde{g}} \check{u}_\e - \frac{ n(n-4)(n^2-4)  \e^4}{\left( \e ^2 + |x|^2
 \right)^{\frac{n+4}{2}}} \right| \leq \begin{cases}
 \frac{O(1)}{\left( \e^2 + r^2 \right)^{\frac{n-4}{2}}} & \hbox{ for } |x| \leq \tilde{\d}, \\
   O(\tilde{\d}^{-3}) & \hbox{ for } \tilde{\d} \leq |x| \leq 2 \tilde{\d}.
 \end{cases}
$$
We can use the latter estimate to control the difference between $u_\e$ and $\check{u}_\e$
by convolving with the Green's function.

\begin{lem}\label{l:diff2}
The following estimate holds, for some constant $C > 0:$
$$
  |\hat{u}_\e - \check{u}_\e| \leq o(1) + C \, \tilde{\d}^{n-3} \min \left\{ |x|^{4-n}, \d^{4-n} \right\}
 = o(1), \ \ \tilde{\d} \rightarrow 0.
$$
\end{lem}

\begin{pf}
By the formula before the lemma we can write that $|\hat{u}_\e - \check{u}_\e| \leq u_1 + u_2$,
where
$$
  u_1(x) = \int_{B_{\tilde{\d}}(0)} G_x(y) \frac{dy}{\left( \e^2 + |y|^2 \right)^{\frac{n-4}{2}}};
  \qquad \quad u_2(x) = \d^{-3} \int_{\tilde{\d} \leq |y| \leq 2 \tilde{\d}} G_x(y) dy.
$$
To estimate $u_1$ we reason as in the proof of Lemma \ref{l:diff}: we divide again into the cases $|x| = O(\e)$
and $|x| \geq C_0 \e$. In the former case we get
$$
  C \int_{|y| \leq \tilde{\d}} \frac{1}{|x-y|^{n-4}} \frac{dy}{\left( \e^2 + |y|^2 \right)^{\frac{n-4}{2}}} =
  C \e^{8-n} \int_{|w| \leq \tilde{\d}/\e} \frac{1}{|\overline{x}-w|^{n-4}} \frac{dw}{\left( 1 + |w|^2 \right)^{\frac{n-4}{2}}},
$$
with $|\overline{x}| = O(1)$. The last integral is  uniformly bounded by $\tilde{\d}^{n-8} \e^{n-8}$, so we
get a quantity of order $o(1)$ as $\tilde{\d} \to 0$.

If the case $|x| \geq C_0 \e$ we write
$$
  \int_{|y| \leq \tilde{\d}} \frac{1}{|x-y|^{n-4}} \frac{dy}{\left( \e^2 + |y|^2 \right)^{\frac{n-4}{2}}}
  = |x|^{8-n} \int_{|y| \leq \frac{\tilde{\d}}{|x|}} \frac{dw}{\left| \frac{x}{|x|} - w \right|^{n-4}
  \left( \frac{\e^2}{|x|^2} + |w|^2 \right)^{\frac{n-4}{2}}} \leq \tilde{\d}^{n-8}.
$$
Therefore we get a uniform bound on $u_1$ of order $o(1)$ as $\tilde{\d} \to 0$.

Turning to $u_2$, one can distinguish the cases $|x| \leq 2 \tilde{\d}$ and $|x| > 2 \tilde{\d}$.
In the former one finds $|u_2(x)| \leq C \tilde{\d}$. In the latter
$$
  |u_2(x)| \leq C \tilde{\d}^{n-3} |x|^{4-n}.
$$
The bounds on $u_1$ and $u_2$ yield the conclusion.
\end{pf}

To estimate the quotient of $\hat{u}_\e$, we have by definition of $\hat{u}_\e$ that
\begin{eqnarray*}
\int_M \hat{u}_\e P_{\tilde{g}} \hat{u}_\e dv_{\tilde{g}} & = & \int_M \hat{u}_\e
  \frac{\eta(x) b_n \e^4}{\left( \e ^2 + |x|^2 \right)^{\frac{n+4}{2}}}  dv_{\tilde{g}} \\
  & = &  \int_M \left( \tilde{\chi}_{\tilde{\d}}(u_\e + \beta) + (1 - \tilde{\chi}_{\tilde{\d}}) \ov{G}_{x_0}
  + (\hat{u}_\e - \check{u}_\e) \right)
    \frac{\eta(x) b_n \e^4}{\left( \e ^2 + |x|^2 \right)^{\frac{n+4}{2}}}  dv_{\tilde{g}}  \\
    & =: & I_1 + I_2 + I_3.
\end{eqnarray*}
We next estimate each of these three terms. Concerning $A_1$ we have
$$
  I_1 = b_n \e^4 \left( \int_M u_\e^{\frac{2n}{n-4}} dv_{\tilde{g}} + \beta \int_M u_\e^{\frac{n+4}{n-4}}
    dv_{\tilde{g}}  + O(1) \right), \ \ \tilde{\d} \to 0.
$$
For $I_2$, since $ (1 - \tilde{\chi}_{\tilde{\d}})$ vanishes in a $\tilde{\d}$-neighborhood
of $p$ we have simply that
$$
  I_2 = \e^4 O(1).
$$
For $I_3$ we can use Lemma \ref{l:diff2} to find that
$$
  |I_3| \leq C o(1) \int_M
      \frac{\eta(x) b_n \e^4}{\left( \e ^2 + |x|^2 \right)^{\frac{n+4}{2}}}  dv_{\tilde{g}}.
$$
Therefore, we obtain
\begin{equation}\label{eq:nuelow}
 \mathcal{N}(\hat{u}_\e) = b_n \e^4 \left( \int_M u_\e^{\frac{2n}{n-4}} dv_{\tilde{g}} +
 \beta (1 + o(1)) \int_M u_\e^{\frac{n+4}{n-4}} dv_{\tilde{g}}  + O(1) \right).
\end{equation}
On the other hand for the denominator we have
$$
 \int_M \hat{u}_\e^{\frac{2n}{n-4}}
   dv_{\tilde{g}} =\int_{B_{\tilde{\d}}(x_0)} (u_\e + \beta)^{\frac{2n}{n-4}}
     dv_{\tilde{g}} + O(1).
$$
Since in $B_{\tilde{\d}}(x_0)$, $\beta$ is bounded by $u_\e$, we have that
\begin{eqnarray*}
  \int_M \hat{u}_\e^{\frac{2n}{n-4}}
     dv_{\tilde{g}} & = & \int_{B_{\tilde{\d}}(x_0)} u_\e^{\frac{2n}{n-4}}
            dv_{\tilde{g}} + \frac{2n}{n-4} \beta \int_{B_{\tilde{\d}}(x_0)} {u}_\e^{\frac{n+4}{n-4}}
       dv_{\tilde{g}} \\ & + & \beta^2 \int_{B_\d(x_0)}  O(u_\e^{\frac{8}{n-4}})
                   dv_{\tilde{g}} + O(1) \\ & = &
                  \int_{M} u_\e^{\frac{2n}{n-4}}
                              dv_{\tilde{g}} + \frac{2n}{n-4} \beta \int_{M} {u}_\e^{\frac{n+4}{n-4}}
                         dv_{\tilde{g}} + \beta^2 \int_{B_{\tilde{\d}}(x_0)}  O(u_\e^{\frac{8}{n-4}})
                                     dv_{\tilde{g}} + O(1).
\end{eqnarray*}
Therefore one finds that
$$
    \mathcal{F}_{\tilde{g}}(\hat{u}_\e) = \frac{b_n \e^4 \left( \int_M u_\e^{\frac{2n}{n-4}} dv_{\tilde{g}} +
     \beta (1 + o_{\tilde{\d}}(1)) \int_M u_\e^{\frac{n+4}{n-4}} dv_{\tilde{g}}  + O_{\tilde{\d}}(1)
     \right)}{\left( \int_{M} u_\e^{\frac{2n}{n-4}}
                                   dv_{\tilde{g}} + \frac{2n}{n-4} \beta \int_{M} {u}_\e^{\frac{n+4}{n-4}}
                              dv_{\tilde{g}} + \beta^2 \int_{B_{\tilde{\d}}(x_0)}  O(u_\e^{\frac{8}{n-4}})
                                          dv_{\tilde{g}} + O_{\tilde{\d}}(1) \right)^{\frac{n-4}{n}}}.
$$
We now notice that the following asymptotics hold
$$
  \int_{M} u_\e^{\frac{2n}{n-4}} dv_{\tilde{g}} \simeq \e^{-n}; \qquad \quad
  \int_{M} {u}_\e^{\frac{n+4}{n-4}}  dv_{\tilde{g}} \simeq \e^{-4}; \qquad \quad
  \int_{B_{\tilde{\d}}(x_0)}  O(u_\e^{\frac{8}{n-4}})  dv_{\tilde{g}} \simeq \e^{n-8}.
$$
These and a Taylor expansion of the denominator in $\mathcal{F}_{\tilde{g}}(\hat{u}_\e)$ imply
$$
  \mathcal{F}_{\tilde{g}}(\hat{u}_\e) = S_n \left( 1 - \beta (1 + o(1))
  \frac{\int_M u_\e^{\frac{n+4}{n-4}} dv_{\tilde{g}}}{\int_{M}
  u_\e^{\frac{2n}{n-4}} dv_{\tilde{g}}}  \right).
$$
This completes the proof of (\ref{Fbarp2}).  The proof of $(i') - (iii')$ is the same as in the proof of Proposition \ref{p:testhighd}.
\end{pf}

%%%%%%%%%%%%%%%%%%%%%%%%%%
\section{Sequential convergence of the flow} \label{converge}
%%%%%%%%%%%%%%%%%%%%%%%%%%%%

In this section we prove the main existence result: under the assumptions of Proposition \ref{p:testhighd} or \ref{p:testlowd}, we show the flow converges (up to choosing a suitable sequence of times)
to a solution of the $Q$-curvature equation.

\begin{thm} \label{flowcon} Let $(M^n,\bar{g})$  be a closed Riemannian manifold of dimension $n \geq 5$ which is not conformally equivalent to the
standard sphere. Suppose that  \\

\noindent $(i)$ $Q_{\bar{g}}$ is semi-positive, \\

\noindent $(ii)$ $R_{\bar{g}} \geq 0$.  \\

%\noindent In addition, assume the dimension $n = 5,6,7$, {\em or} that $n \geq 8$ and $(M^n,\bar{g})$ is locally conformally flat.  \\

Let $g_0 = h$, where $h$ is the metric constructed in Proposition \ref{p:testlowd} (when $5 \leq n \leq 7$, or $\bar{g}$ is locally conformally flat and $n \geq 5$) or Proposition \ref{p:testhighd} (when $n \geq 8$ and $\bar{g}$ is not locally conformally flat).  Then the flow (\ref{uflow}) has a solution for all time satisfying
\begin{align} \label{belowL2}
\int u^2\ dv_0 \geq C_0
\end{align}
for some constant $C_0 > 0$.  Moreover, it is possible to choose s sequence of times $t_j \nearrow \infty$ such that $u_j = u_j(t_j, \cdot)$ converges weakly in $W^{2,2}(M^n)$ to a
smooth solution $u > 0$ of
\begin{align} \label{Endgame}
P_{g_0} u = \bar{\mu} \, u^{\frac{n+4}{n-4}},
\end{align}
where $\bar{\mu} > 0$.  In particular, $g_{\infty} = u^{\frac{4}{n-4}}g_0$ defines a metric with positive scalar curvature and constant positive $Q$-curvature.
\end{thm}

\begin{pf}
If we take our initial metric $g_0$ to be the metric in the conclusion of Proposition \ref{p:testhighd} or \ref{p:testlowd}, then by Proposition \ref{LTE} we know the flow (\ref{uflow}) exists for all time.  In addition, by same Propositions we know
\begin{align} \label{beginsmall}
\mathcal{F}_{g_0}[ u_0 ] \leq S_n - \epsilon_0,
\end{align}
where $u_0 \equiv 1$ is our initial datum for the flow and $\epsilon_0 > 0$.  It follows from Lemma \ref{VLem} that
\begin{align} \label{smalltime}
\mathcal{F}_{g_0}[u] = \frac{ \int u (P_{g_0} u)\ dv_0 }{ \Big( \int u^{\frac{2n}{n-4}}\ dv_0 \Big)^{\frac{n-4}{n}} } \leq S_n - \epsilon_0
\end{align}
for all times.

Recall the Euclidean Paneitz-Sobolev constant is
\begin{align*}
S_n = \inf_{\varphi \in C_0^{\infty}(\mathbb{R}^n)} \frac{ \int (\Delta_0 \varphi)^2\ dx }{ \Big( \int |\varphi|^{\frac{2n}{n-4}}\ dx \Big)^{\frac{n-4}{n}}}.
\end{align*}
On the compact Riemannian manifold $(M,g_0)$, given $\delta > 0$ we can use a cut-and-paste argument to prove that
\begin{align*}
\Big( \int |\varphi|^{\frac{2n}{n-4}}\ dv_0 \Big)^{\frac{n-4}{n}} \leq \big( S_n^{-1} + \delta \big) \int (\Delta_{g_0} \varphi)^2\ dv_0 + C_{\delta} \int \varphi^2\ dv_0,
\end{align*}
which implies
\begin{align} \label{Sob}
\Big( \int |\varphi|^{\frac{2n}{n-4}}\ dv_0 \Big)^{\frac{n-4}{n}} \leq \big( S_n^{-1} + 2\delta \big) \int \varphi (P_{g_0} \varphi)\ dv_0 + C_{\delta}' \int \varphi^2\ dv_0.
\end{align}
Plugging (\ref{smalltime}) into the Sobolev inequality (\ref{Sob}) gives
\begin{align} \label{low1} \begin{split}
 \Big( \int u^{\frac{2n}{n-4}}\ dv_0 \Big)^{\frac{n-4}{n}} &\leq \big( S_n^{-1} + 2\delta \big) \int u (P_{g_0} u)\ dv_0 + C_{\delta}' \int u^2\ dv_0 \\
 &\leq \big( S_n^{-1} + 2\delta \big)\big( S_n - \epsilon_0 \big) \Big( \int u^{\frac{2n}{n-4}}\ dv_0 \Big)^{\frac{n-4}{n}} + C_{\delta}' \int u^2\ dv_0.
 \end{split}
\end{align}
If we take $\delta = \epsilon_0 /10$, then the first term on the right-hand side can be absorbed into the left-hand side, and we get
\begin{align} \label{low2}
\Big( \int u^{\frac{2n}{n-4}}\ dv_0 \Big)^{\frac{n-4}{n}} \leq C(\epsilon_0) \int u^2\ dv_0.
\end{align}
Since the l.h.s. is just a power of the conformal volume (which is non-decreasing), we conclude
\begin{align*}
\int u^2\ dv_0 \geq C_0 > 0
\end{align*}
for all time, as claimed.

By Lemma \ref{VLem} and Corollary \ref{STCor} we can choose a sequence of times $t_j \nearrow \infty$ such that
$u_j = u(t_j, \cdot)$ and $\mu_j = \mu(t_j)$
satisfy
\begin{align} \label{utimesj} \begin{split}
\mu_j &\nearrow \bar{\mu}, \\
u_j &\rightharpoonup u \mbox{ weakly in } W^{2,2}(M^n), \\
u_j &\rightarrow u \mbox{ strongly in } L^2(M^n), \\
f_j &= - u_j + \mu_j P_{g_0}^{-1} (u_j^{\frac{n+4}{n-4}} ) \rightarrow 0 \mbox{ strongly in } W^{2,2}(M^n).
\end{split}
\end{align}
It follows that $u \geq 0$ satisfies
\begin{align} \label{limitu}
u = \bar{\mu} P_{g_0}^{-1} \big( u^{\frac{n+4}{n-4}} \big),
\end{align}
and by elliptic regularity $u$ is a strong solution of
\begin{align*}
P_{g_0} u = \bar{\mu} u^{\frac{n+4}{n-4}}.
\end{align*}
By the strong maximum principle Theorem \ref{SMP}, in fact $u > 0$.  This completes the proof of the theorem.
\end{pf}

\bibliography{Qflow_references}

\end{document}